\documentclass[2p]{article}
\usepackage{amssymb,amsmath,amsthm}
\usepackage{indentfirst}
\usepackage{exscale}
\usepackage{relsize}
\usepackage[numbers,sort&compress]{natbib}

\usepackage{geometry}
\newcommand{\R}{\mathbb{R}}
\newtheorem{theorem}{Theorem}[section]

\theoremstyle{definition}
\newtheorem{definition}[theorem]{Definition}

\newtheorem{lemma}[theorem]{Lemma}

\newtheorem{corollary}[theorem]{Corollary}

\allowdisplaybreaks
\theoremstyle{remark}
\newtheorem{remark}[theorem]{Remark}
\numberwithin{equation}{section}

\begin{document}

\title{\Large\bf{Multiplicity of solutions for a class of  nonhomogeneous quasilinear  elliptic system with  locally symmetric condition in $\mathbb{R}^N$}
 }
\date{}
\author {\ Cuiling Liu$^{1}$, \ Xingyong Zhang$^{1,2}$\footnote{Corresponding author, E-mail address: zhangxingyong1@163.com}, \ Liben Wang$^{3}$ \\
      {\footnotesize $^{1}$Faculty of Science, Kunming University of Science and Technology, Kunming, Yunnan, 650500, P.R. China.}\\
      {\footnotesize $^{1,2}$Research Center for Mathematics and Interdisciplinary Sciences, Kunming University of Science and Technology,}\\
 {\footnotesize Kunming, Yunnan, 650500, P.R. China.}\\
 {\footnotesize $^{3}$School of Computer Science and Technology, Dongguan University of Technology, Dongguan, Guangdong, 523808, P.R. China.}\\
 }
 \date{}
 \maketitle

 \begin{center}
 \begin{minipage}{15cm}

\small  {\bf Abstract:}
This paper is concerned with a class of nonhomogeneous quasilinear  elliptic system driven by the locally symmetric potential and the small continuous perturbations in the whole-space $\mathbb{R}^N$. By a variant of Clark's theorem without the global symmetric condition and a Moser's iteration technique, we obtain the existence of  multiple solutions when the nonlinear term  satisfies some growth conditions only in a circle with center 0 and the perturbation term is any continuous function with a small parameter and no any growth hypothesis.

\end{minipage}
 \end{center}
  \allowdisplaybreaks
 \vskip2mm
\section{Introduction}\label{section 1}
In this paper, we are interested in  the following  non-homogeneous  elliptic problem with  perturbation:
 \begin{equation}\label{eq1}
 \left\{
  \begin{array}{ll}
 -\Delta_{\Phi_1}u+V_{1}(x)\phi_{1}(|u|)u=F_u(x,u,v)+\varepsilon k(x)G_u(u,v), &x\in \mathbb{R}^N,\\
 -\Delta_{\Phi_2}v+V_{2}(x)\phi_{2}(|v|)v=F_v(x,u,v)+\varepsilon k(x)G_v(u,v), &x\in \mathbb{R}^N,\\
  u\in W^{1,\Phi_{1}}(\mathbb{R}^N),v\in W^{1,\Phi_{2}}(\mathbb{R}^N),
    \end{array}
 \right.
 \end{equation}
where $N>2$,  $\Delta_{\Phi_i}u=\mbox{div}(\phi_i(|\nabla u|)\nabla u)$, $i=1, 2$, $\Phi_i(t):=\int_{0}^{|t|}s\phi_i(s)ds, \ t\in \mathbb{R}$, and $\phi_i  :(0,+\infty)\rightarrow(0,+\infty)$ are two functions
satisfying
\begin{itemize}
\item[$(\phi_1)$]
$\phi_i\in C^1(0,+\infty)$, $t\phi_i(t)\rightarrow0$ as
$t\rightarrow0$, $t\phi_i(t)\rightarrow+\infty$ as
$t\rightarrow+\infty$;
\item[$(\phi_2)$]
$t\rightarrow t\phi_i(t)$ are strictly increasing;
\item[$(\phi_3)$]
$1<l_i:=\inf_{t>0}\frac{t^2\phi_i(t)}{\Phi_i(t)}\leq\sup_{t>0}\frac{t^2\phi_i(t)}{\Phi_i(t)}
=:m_i<\min\{N,l_{i}^{\ast},l_{i}\min\{l_{1},l_{2}\}\}$,
where $l_{i}^{\ast}=\frac{l_{i}N}{N-l_{i}}$ and $l_{i}^{\ast}m_i(l_i+1)-l_{i}^{2}(l_{i}^{\ast}+m_i)\geq 0$;
\item[$(\phi_4)$]
there exist positive constants $C_{i,1}$ and $C_{i,2}$, $i=1,2$, such that
$$
 C_{i,1}|t|^{l_i}\le \Phi_i(t)\le C_{i,2}|t|^{l_i},\ \ \forall \;|t|<1.
$$
\end{itemize}
Quasilinear elliptic system of the form \eqref{eq1} includes the following class of quasilinear elliptic scalar  equation:
\begin{equation}\label{eq2}
\left\{
\begin{array}{ll}
-\mbox{div}(\phi(|\nabla{u}|)\nabla{u})+V(x)\phi(|u|)u=f(x,u)+\varepsilon k(x)g(u), &x\in \mathbb{R}^N,\\
u\in W^{1,\Phi}(\mathbb{R}^N),
\end{array}
\right.
\end{equation}
which corresponds to the special case $\phi_2=\phi_1=:\phi$, $v=u$, $V_{2}(x)=V_{1}(x):=V(x)$ and $F(x,u,v)=F(x,v,u)$.   The problems such as \eqref{eq1} or \eqref{eq2} often appear  in the model of nonlinear elasticity, plasticity and  non-Newtonian fluid (see \cite{Fuchs2000,Ruzicka2000,Fukagai1995,Zhikov1987}  in details).
\par
When $\varepsilon=0$, there are fruitful results concerning the existence and multiplicity of solutions for system \eqref{eq1} or equations like  \eqref{eq2} on a bounded domain $\Omega\subset \R^N$. We refer readers to \cite{Clement2000,wang2017,wang2017-3,Huidobro1999} and the references therein. For quasilinear problem \eqref{eq1} on the whole space $\mathbb{R}^N$,  because of the lack of compactness of the Sobolev embedding,  there are limited  works on the existence and multiplicity of solutions for system \eqref{eq1}, for example, \cite{wang2017-2,Zhang2021,Wang2016}. In \cite{wang2017-2}, Wang-Zhang-Fang considered the existence of ground state solution for system \eqref{eq1} with $\varepsilon=0$. Under some called $(\phi_{1},\phi_{2})$-superlinear Orlicz-Sobolev conditions, a  ground state solution for system  \eqref{eq1} is obtained via  a variant mountain pass lemma.  In \cite{Wang2016}, Wang-Zhang-Fang considered the existence and multiplicity of solutions for  system \eqref{eq1} with $\varepsilon=0$. One nontrivial weak solution and infinitely many solutions of system \eqref{eq1} are obtained via the least action principle and the genus theory, respectively, where $F$ satisfies a sub-linear Orlicz-Sobolev growth and global symmetric condition. In \cite{Zhang2021}, Zhang-Liu  considered the existence of solutions for system \eqref{eq1} with a parameter and $\varepsilon=0$. Under some growth assumptions only in a circle for $F$,  a nontrivial solution of system \eqref{eq1} is obtained by the mountain pass lemma and the Moser's iteration technique.
\par
If $\varepsilon$ is allowed to be nonzero, in \cite{Kajikiya2013}, Kajikiya investigated the existence of multiple solutions for the second order elliptic equation with perturbation, which is a special example of equation \eqref{eq2} with $\phi(t)=1$ and $V(x)=0$,
\begin{equation}\label{eq3}
 \left\{
  \begin{array}{ll}
 -\Delta u=f(x,u)+\varepsilon g(x,u), &x\in \Omega,\\
 u=0, &x\in\partial\Omega,
 \end{array}
 \right.
\end{equation}
where $\Omega$ is a bounded domain in $\mathbb{R}^N$ with smooth boundary $\partial \Omega$ and $\varepsilon$ is a small parameter. Kajikiya assumed that the nonlinear term $f(x,u)$ is sublinear and odd only near $u=0$ and $g(x,u)$ is any continuous function near  $u=0$ without any growth hypothesis. Such problem has two difficulties: (i) the variational functional is not well defined and continuously differentiable on the working space because all of the assumptions of the nonlinear terms $f$ and $g$ are local, just near the origin; (ii) the symmetry of equation is broken by the perturbation term $\varepsilon g$. These two difficulties cause that the classical Clark's Theorem (\cite{Rabinowitz1986}, Theorem 9.1) can not work well for equation (\ref{eq3}). To overcome the difficulties,  Kajikiya developed  a new version of Clark's theorem and then together with cut-off technique,  obtained that equation \eqref{eq3} has at least $k$ distinct solutions for any given $k\in \mathbb N$.
\par
Motivated by \cite{Wang2016,Kajikiya2013}, the main aim of this paper is to extend the result in \cite{Kajikiya2013} to  the quasilinear system \eqref{eq1}. We assume that  $F$ is sub-$\min\{l_1,l_2\}$ growth and locally even with respect to $(u,v)$ and  the perturbation term $G$ is any continuous function with a small coefficient $\varepsilon$. We obtain a multiplicity result  of nontrivial solutions for  system \eqref{eq1}.
Our results develop those results in  \cite{Wang2016,Kajikiya2013} from the following two aspects:
\begin{itemize}
\item[(\uppercase\expandafter{$I_{1}$})]
Different from  \cite{Kajikiya2013}, we work in the whole-space $\mathbb{R}^N$ rather than in a bounded domain $\Omega\subset\mathbb{R}^N$ and we consider the more general quasilinear $(\phi_1,\phi_2)$-Laplacian system  \eqref{eq1}. The non-homogeneity of $\Phi_i,i=1,2$ and the mutual coupling of $u$ and $v$ make our proofs become  much difficult, especially, in the proofs that {\rm(P-S)} condition holds for functional $I_{\varepsilon}$ (Lemma 3.4-Lemma 3.6 below) and the proofs of relation between $\|u\|_{\infty}$, $\|v\|_{\infty}$ and $\|u\|_{l_{1}^{\ast}}$, $\|v\|_{l_{2}^{\ast}}$ (Lemma 3.9 below) by Moser's iteration technique. We overcome these difficulties by using some properties of $\Phi_i,i=1,2$ sufficiently and scaling some inequalities carefully.
\item[(\uppercase\expandafter{$I_{2}$})]
We consider the case that $\varepsilon\not=0$ and the continuity is an unique assumption for the perturbation term $G$, which implies that the local symmetry of system \eqref{eq1} is broken if $G(u,v)$ is not even. This case is not considered in \cite{Wang2016} and  our conditions are still different from those in \cite{Wang2016} even if $\varepsilon=0$, because all of our assumptions on $F(x,u,v)$ and $G(u,v)$ are local,  just near the origin, but those  assumptions except for $(F_2)$ in \cite{Wang2016} are global.
\end{itemize}
\par
To state our result, we introduce the following  assumptions concerning $F$, $V_{i}$, $K$ and $G$:
\begin{itemize}
\item[$(F_0)$]
$F: \mathbb{R}^N\times [-\sqrt{2}\delta, \sqrt{2}\delta] \times [-\sqrt{2}\delta, \sqrt{2}\delta]\rightarrow \mathbb{R}$ is a $C^1$ function for some $\delta>0$,  such that
$F(x,0,0)=0$ for all $x\in \mathbb{R}^N$ and $F(x,u,v)>0$ for all $0<|(u,v)|<\sqrt{2}\delta$ and $x\in  \mathbb{R}^{N}$;
\item[$(F_1)$]
there exist two  constants $C_1, C_2>0$ such that
\begin{equation*}
 \left\{
  \begin{array}{l}
 |F_u(x,u,v)|\leq C_1\left(|u|^{l_{1}r-1}+|v|^{\frac{l_{2}(l_{1}r-1)}{l_{1}}}\right),\\
 |F_v(x,u,v)|\leq C_2\left(|v|^{l_{2}r-1}+|u|^{\frac{l_{1}(l_{2}r-1)}{l_{2}}}\right)
    \end{array}
 \right.
 \end{equation*}
for all $|(u,v)|<\delta$  and $x\in  \mathbb{R}^{N}$, where
$r\in \left(\max\left\{\frac{1}{l_{1}},\frac{1}{l_{2}}\right\}, \min\left\{\frac{l_{1}}{m_{1}},\frac{l_{2}}{m_{2}}\right\}\right)$;
\item[$(F_2)$]
there exists  a constant $\alpha\in \left[\max\{l_{1},l_{2}\}r,\min\{l_{1},l_{2}\}\right)$ such that
{\small\begin{eqnarray*}
      F_u(x,u,v)u+F_v(x,u,v)v-\alpha F(x,u,v)
\leq a_{1}(x)|u|^{l_{1}}
    +a_{2}(x)|v|^{l_{2}}
\end{eqnarray*}}
for all $|(u,v)|<\delta$  and $x\in  \mathbb{R}^{N}$, where $a_{1}(x)=\frac{\min\left\{C_{1,1},\Phi_{1}(1)\right\}\left(l_{1}-\alpha\right)}{l_{1}}V_{1}(x)$
and
$a_{2}(x)=\frac{\min\left\{C_{2,1},\Phi_{2}(1)\right\}\left(l_{2}-\alpha\right)}{l_{2}} V_{2}(x)$;
\item[$(F_3)$]
$\lim_{|(u,v)|\rightarrow 0}\left(\inf_{x \in \mathbb{R}^{N}}\frac{F(x,u,v)}{|u|^{\alpha}+|v|^{\alpha}}\right)=+\infty$;
\item[$(F_4)$]
 $F(x,-u,-v)=F(x,u,v)$
for all $|(u,v)|<\delta$  and $x\in  \mathbb{R}^{N}$;
\end{itemize}
\begin{itemize}
\item[$(V)$]
$V_{i} \in C(\mathbb{R}^{N},\mathbb{R})$,
$V_{i}(x)^{-1}\in L^{\frac{r}{1-r}}(\mathbb{R}^{N})$
\mbox{and} there exists a constant $C_{3}$ such that
$0<C_{3}<V_{i}(x)$,\; $i=1,2$;
\end{itemize}
\begin{itemize}
\item[$(K)$]
$k(x)\in L^{1}(\mathbb{R}^N)\cap L^{\infty}(\mathbb{R}^N)$.
\end{itemize}
\begin{itemize}
\item[$(G)$]
$G: [-\sqrt{2}\delta, \sqrt{2}\delta] \times [-\sqrt{2}\delta, \sqrt{2}\delta]\rightarrow \mathbb{R}$ is a $C^1$ function such that
                $G(0,0)=0$.
\end{itemize}
\par
Our main result is as follows.
\par
\begin{theorem}\label{theorem1.1}
  Assume that $(\phi_1)$--$(\phi_4)$, $(F_0)$--$(F_4)$, $(V)$, $(G)$ and $(K)$ hold.
 Then
 for any $ k \in \mathbb{N}$ and  any $b_{1}>0$, $b_{2}>0$, there exists an  $\sigma_{0}(k,b_{1},b_{2})>0$ such that when $|\varepsilon|\leq \sigma_{0}(k,b_{1},b_{2})$,  system \eqref{eq1} possesses at least $k$ distinct solutions whose
 $L^{\infty}$-norms are less than $\frac{\sqrt{2}\delta}{2}$.
\end{theorem}
\par
By Theorem 1.1, it is easy to obtain the following corollary.
\begin{corollary}\label{corollay1.2}
{\it Assume that $(\phi_1)$--$(\phi_4)$, $(F_0)$--$(F_4)$ and $(V)$ hold. Then system \eqref{eq1} with $\varepsilon=0$ possesses infinitely many distinct solutions whose
 $L^{\infty}$-norms are less than $\frac{\sqrt{2}\delta}{2}$.}
\end{corollary}

\begin{remark}\label{remark1.2}
In Theorem 1.1, $F$ does not need to satisfy any growth condition at infinity. So $F$ is allowed to be subcritical, critical or supercritical, for example,
\begin{itemize}
\item[$(1)$] the subcritical case:
    $F(x,t,s)=a(x)(|t|^{r_{1}}+|s|^{r_{2}}+|t|^{r_{3}}|s|^{r_{4}}+|t|^{r_{5}}+|s|^{r_{6}})$,
where $l_{1}r\leq r_{1}<\alpha<l_{1}$, $l_{1}r<\alpha$, $l_{2}r\leq r_{2}<\alpha<l_{2}$, $l_{2}r<\alpha$,
$r_{3}\geq l_{1}$, $r_{4}\geq l_{2}$, $l_{1}\leq r_{5}< l_{1}^{\ast}$, $l_{2}\leq r_{6}< l_{2}^{\ast}$, $a(x)$ is continuous and $0<a_{0}\leq a(x)\leq a_{1}$;
\item[$(3)$] the critical case:
 $F(x,t,s)=a(x)(|t|^{r_{1}}+|s|^{r_{2}}+|t|^{l_{1}^{\ast}}+|s|^{l_{2}^{\ast}})$,
where $l_{1}r\leq r_{1}<\alpha<l_{1}$, $l_{1}r<\alpha$, $l_{2}r\leq r_{2}<\alpha<l_{2}$, $l_{2}r<\alpha$,
$a(x)$ is continuous and $0<a_{0}\leq a(x)\leq a_{1}$;
\item[$(2)$] the supercritical case:
$F(x,t,s)=a(x)(|t|^{r_{1}}+|s|^{r_{2}}+|t|^{r_{7}}+|s|^{r_{8}})$,
where $l_{1}r\leq r_{1}<\alpha<l_{1}$, $l_{1}r<\alpha$, $l_{2}r\leq r_{2}<\alpha<l_{2}$, $l_{2}r<\alpha$,
$r_{7}> l_{1}^{\ast}$, $r_{8}> l_{2}^{\ast}$, $a(x)$ is continuous and $0<a_{0}\leq a(x)\leq a_{1}$.
\end{itemize}
\end{remark}

 \section{Preliminaries}\label{section 2}
This section focuses on a survey of concepts and results from Orlicz and Orlicz-Sobolev spaces that
will be used in the text  and introduce the result from variational methods.  For a deeper understanding of these concepts, we refer readers for more details to the books
\cite{Adams2003,Rabinowitz1986, M. M. Rao2002, M.A.Krasnoselski1961}.
\par
\noindent
\begin{definition}\cite{Adams2003}\label{definition2.1}
The function defined on $\mathbb{R}$ by $B(t)=\int_{0}^{|t|}b(s)ds$ is called an $N$-function, if $b$ : $[0,+\infty)\rightarrow [0,+\infty)$ be a right continuous, monotone increasing function with
\begin{itemize}
\item[$\rm(1)$] $b(0)=0$;
\item[$\rm(2)$] $\lim_{t\rightarrow +\infty} b(t)=+\infty$;
\item[$\rm(3)$] $b(t)>0 $ whenever $t>0$.
\end{itemize}
\end{definition}
\par
\noindent
\begin{definition}\cite{Adams2003}\label{definition2.1.1}
The $N$-function $B$ satisfies a $\Delta_2$-condition globally (or near infinity) if
$$
  \sup_{t>0}\frac{B(2t)}{B(t)}<+\infty \ \ \left(\mbox{ or } \limsup_{t\rightarrow \infty}\frac{B(2t)}{B(t)}<+\infty\right),
$$
which implies that there exists a constant $K>0$, such that $B(2t)\leq K B(t)$ for all $t\geq0$ (or $t\geq t_0>0$).
Moreover, $B$ satisfies a $\Delta_2$-condition globally (or near infinity) if and only if for any $c\geq 1$, there exists a constant $K_c>0$ such that $B(ct)\leq K_c B(t)$ for all $t\geq0$ (or $t\geq t_0>0$).
\end{definition}
\par
\noindent
\begin{definition}\cite{Adams2003}\label{definition2.2}
For an $N$-function $B$,  define
$$ \widetilde{B}(t)=\int_{0}^{|t|}b^{-1}(s)ds, \quad t\in \mathbb{R},$$
 where $b^{-1}$ is the right inverse of the right derivative $b$ of $B$. Then $\widetilde{B}$ is an $N$-function
called as the complement of $B$.
\end{definition}
\par
\noindent
\begin{lemma}\label{lemma2.6.1}
$N$-function $B$ and $\widetilde{B}$, which is the complement of $B$, satisfy the Young's inequality (see \cite{Adams2003,M. M. Rao2002})
\begin{equation}\label{2.1.1}
st\leq B (s)+\widetilde{B}(t), \quad s, t\geq 0
\end{equation}
and the following inequality (see \cite[Lemma A.2]{Fukagai2006})
 \begin{equation}\label{2.1.2}
 \widetilde{B}(b(t))\leq B(2t), \quad t\geq 0.
 \end{equation}
\end{lemma}
\par
Now, we recall the Orlicz space $L^{B}(\Omega)$ associated with $N$-function $B$.
The Orlicz space $L^{B}(\Omega)$ is the vectorial space of the measurable functions $u: \Omega\rightarrow \mathbb{R}$ satisfying
$$\int_{\Omega}B(|u|)dx<+\infty,$$
where $\Omega \subset \mathbb{R}^N$ is an open set. It was obvious that $L^{B}(\Omega)$ is a Banach space endowed with Luxemburg norm
$$
\|u\|_{B}:=\inf \left\{\lambda >0: \int_{\Omega}B \left(\frac{u}{\lambda}\right)dx\le 1\right\}.
$$
The corresponding Orlicz-Sobolev space (see \cite{Adams2003,M. M. Rao2002}) is defined by
$$
W^{1, B}(\Omega):=\left\{u \in L^{B}(\Omega): \frac{\partial u}{\partial x_i} \in L^{B}(\Omega), i=1,\cdots, N\right\}
$$
with the norm
\begin{eqnarray}\label{xxxxxx1}
\|u\|:=\|u\|_{B}+\|\nabla u \|_{B}.
\end{eqnarray}
\par
When $\Omega=\mathbb{R}^N$, the Orlicz-Sobolev space $W^{1, B}(\mathbb{R}^N)$ is the completion of $C_{0}^{\infty}(\mathbb{R}^N)$ under the norm \eqref{xxxxxx1}.
Moreover,
if we assume that $V \in C(\mathbb{R}^{N},\mathbb{R})$ which satisfies
$V_{0}:=\inf_{x\in \mathbb{R}^{N}}V(x)>0$.
Then, in order to deal with the system \eqref{eq1}, we introduce the subspace $X$ of $W^{1,B}(\mathbb{R}^N)$, which is  defined by
\begin{eqnarray}\label{dddc1}
  X=\left\{u\in W^{1,B}(\mathbb{R}^N)\Big| \int_{\mathbb{R}^N}V(x)B(|u|)dx<\infty\right\}
\end{eqnarray}
with the norm
\begin{eqnarray}\label{2.1.4}
  \|u\|_{1,B}=\|\nabla u\|_B+\|u\|_{B,V},
\end{eqnarray}
where
$$
  \|u\|_{B,V}=\inf\left\{\alpha>0\Big|\int_{\mathbb{R}^N} V(x)B\left(\frac{|u|}{\alpha}\right)dx\le 1\right\}.
$$
It is easy to see that $(X,\|\cdot\|)$ is a separable and reflexive Banach space (see \cite{Liu-Shibo2019}).
\par
In what follows,  we recall some inequalities and lemmas which we will use. For more details, we refer the reader to the references \cite{Huang2022,Adams2003,Fukagai2006,Liu-Shibo2019}.
\par
\noindent
\begin{lemma}\cite{Adams2003,Fukagai2006}\label{lemma2.6}
For the $N$-function $B$, the following conditions are equivalent:
\begin{itemize}
	\item[$\rm(1)$]
\begin{equation}\label{2.1.5}
1\leq l=\inf_{t>0}\frac{tb(t)}{B(t)}\leq\sup_{t>0}\frac{tb(t)}{B(t)}=m<+\infty;
\end{equation}
\item[$\rm(2)$] let $\zeta_0(t)=\min\{t^l, t^m\}$ and $\zeta_1(t)=\max\{t^l, t^m\}$, $t\geq0$. $B$ satisfies
$$\zeta_0(t)B(\rho)\leq B(\rho t)\leq \zeta_1(t)B(\rho), \quad \forall \rho, t\geq 0;$$
\item[$\rm(3)$] $B$ satisfies a $\Delta_2$-condition globally.
\end{itemize}
\end{lemma}
\noindent
\begin{lemma}\cite{Fukagai2006,Liu-Shibo2019}\label{lemma2.7}
 If $B$ is an $N$-function and \eqref{2.1.5} holds, then $B$ satisfies
$$\zeta_0(\|u\|_{B})\leq\int_{\Omega}B(u)dx\leq\zeta_1(\|u\|_{B}), \quad \forall u \in L^{B}(\Omega).$$
Moreover, if $l>1$, $V(x)\in C(\mathbb{R}^{N})$,
 $V_{0}=\inf_{\mathbb{R}^{N}}V(x)>0$. Then for all $u\in X$,
$$\zeta_0(\|u\|_{B,V})\leq\int_{\mathbb{R}^{N}}V(x)B(u)dx\leq\zeta_1(\|u\|_{B,V}).$$
\end{lemma}
\noindent
\begin{lemma}\cite{Fukagai2006}\label{lemma2.8}
 If $B$ is an $N$-function and \eqref{2.1.5} holds with $l>1$. Let $\widetilde{B}$ be the complement of $B$ and $\zeta_2(t)=\min\{t^{\widetilde{l}},t^{\widetilde{m}}\}$, $\zeta_3(t)=\max\{t^{\widetilde{l}},t^{\widetilde{m}}\}$ for $t\geq0$, where $\widetilde{l}:=\frac{l}{l-1}$ and $\widetilde{m}:=\frac{m}{m-1}$. Then $\widetilde{B}$ satisfies
\begin{itemize}
	\item[$\rm(1)$]
$$\widetilde{m}=\inf_{t>0}\frac{t\widetilde{B}^{'}(t)}{\widetilde{B}(t)}\leq\sup_{t>0}\frac{t\widetilde{B}^{'}(t)}{\widetilde{B}(t)}=\widetilde{l};$$
\item[$\rm(2)$]
$$\zeta_2(t)\widetilde{B}(\rho)\leq \widetilde{B}(\rho t)\leq \zeta_3(t)\widetilde{B}(\rho), \quad \forall  \rho, t\geq 0;$$
\item[$\rm(3)$]
$$\zeta_2(\|u\|_{\widetilde{B}})\leq\int_{\Omega}\widetilde{B}(u)dx\leq\zeta_3(\|u\|_{\widetilde{B}}), \quad \forall u \in L^{\widetilde{B}}(\Omega).$$
\end{itemize}
\end{lemma}
\par
\noindent
\begin{remark}\label{remark 2.9}
By Lemma 2.3 and Lemma 2.5, assumptions $(\phi_1)$--$(\phi_3)$ show that $\Phi_i ~(i=1, 2)$ and $\widetilde{\Phi}_i ~(i=1, 2)$ are $N$-functions satisfying $\Delta_2$-condition globally. Thus $L^{\Phi_i}(\mathbb{R}^N) (i=1, 2)$ and $W^{1, \Phi_i}(\mathbb{R}^N) (i=1, 2)$ are separable and reflexive Banach spaces (see \cite{Adams2003,M. M. Rao2002}).
\vskip2mm
\par
If the  $N$-function $B$ satisfies $\Delta_2$-condition globally, then we can see that
\begin{equation}\label{2.1.3+}
 u_n\rightarrow u \mbox{ in } L^{B}(\Omega) \Longleftrightarrow \int_{\Omega}B(u_n-u)dx\rightarrow 0.
 \end{equation}
 Moreover, a generalized type of H\"{o}lder's inequality (see \cite{Adams2003,M. M. Rao2002})
$$\left| \int_{\Omega}uvdx \right|\leq 2\|u\|_{B}\|v\|_{\widetilde{B}}, \quad \mbox{ for all } u \in L^{B}(\Omega)\mbox{ and } v \in L^{\widetilde{B}}(\Omega)$$
can be gained by applying Young's inequality \eqref{2.1.1}.
\par
Let $N$-function $A$ satisfies $\Delta_2$-condition. Assume that
 \begin{equation}\label{2.1.6}
\overline{\lim_{t\rightarrow 0}}\frac{A(t)}{B(t)}< +\infty\;\;\;\mbox{and}\;\;\;
\mathop{\overline{\lim}}_{|t|\rightarrow +\infty}\frac{A(t)}{B_{\ast}(t)}< +\infty.
\end{equation}
Then the embedding $W^{1,B}(\mathbb{R}^{N})\hookrightarrow L^{A}(\mathbb{R}^{N})$ is continuous.
Moreover, if we assume that
 \begin{equation}\label{2.1.7}
\lim_{|t|\rightarrow 0}\frac{A(t)}{B(t)}< +\infty\;\;\;\mbox{and}\;\;\;
\lim_{|t|\rightarrow +\infty}\frac{A(t)}{B_{\ast}(t)}=0,
 \end{equation}
then the embedding $W^{1,B}(\mathbb{R}^{N})\hookrightarrow L_{loc}^{A}(\mathbb{R}^{N})$ is compact and we call that such $A$ satisfies the subcritical condition (see \cite{Liu-Shibo2019}).
\end{remark}
\vskip2mm
\par
Finally, we   recall a  version of Clark's lemma which
was introduced in \cite{Kajikiya2013}.
\par
\noindent
\begin{lemma}\cite{Kajikiya2013}\label{lemma2.10}
Let $W$ be an infinite dimensional Banach space. For any $\varepsilon\in[0,1]$,  $I_{\varepsilon}\in C(W, \mathbb{R})$. Suppose that  $I_{\varepsilon}(u)$ has a continuous partial derivative $I_{\varepsilon}'$ with respect to $u$ and  satisfies $(A_{1})-(A_{5})$ below.
\begin{itemize}
\item[$(A_{1})$:]
$\inf \{I_{\varepsilon}(u):\varepsilon \in [0,1], u\in W \}>-\infty$;
\item[$(A_{2})$:]
For $u\in W$, $|I_{\varepsilon}(u)-I_{0}(u)|\leq \psi(\varepsilon)$, where $\psi\in C([0,1],\mathbb{R})$ and $\psi(0)=0$, $I_{0}(u)=:I_{\varepsilon}(u)|_{\varepsilon=0}$;
\item[$(A_{3})$:]
$I_{\varepsilon}(u)$ satisfies the (P-S) conditions uniformly on $\varepsilon$, i.e. if a sequence $(\varepsilon_{k},u_{k})$ in $[0,1]\times W$ satisfies that
$\sup_{k}|I_{\varepsilon_{k}}(u_{k})|<\infty$ and $I_{\varepsilon_{k}}'(u_{k})$ converges to zero, then
$(\varepsilon_{k},u_{k})$ has a convergent subsequence;
\item[$(A_{4})$:]
$I_{0}(u)=I_{0}(-u)$ for $u\in W$ and $I_{0}(0)=0$;
\item[$(A_{5})$:]
For $u\in W\backslash\{0\}$, there exists a unique $t(u)>0$ such that $I_{0}(tu)<0$ if $0<|t|<t(u)$ and $I_{0}(tu)\geq 0$ if $|t|\geq t(u)$.
\end{itemize}
 Denote
$$
S^{k}:=\{x\in \mathbb{R}^{k+1}:|x|=1\},
$$
$$
A_{k}:=\{\alpha\in C(S^{k},W):\alpha \;\;\mbox{is odd}\},
$$
$$
d_{k}:=\inf_{\alpha\in A_{k}}\max_{x\in S^{k}}I_{0}(\alpha(x)).
$$
Let $k\in \mathbb{N}\backslash \{0\}$ satisfying $d_{k}<d_{k+1}$. Then there exist two constants $\varepsilon_{k+1}$, $c_{k+1}$ such that $0<\varepsilon_{k+1}\leq 1$, $d_{k+1}\leq c_{k+1}<-\psi(\varepsilon)$ for $\varepsilon \in [0,\varepsilon_{k+1}]$  and for any $\varepsilon \in [0,\varepsilon_{k+1}]$, $I_{\varepsilon}(\cdot)$ has a critical value in the interval $[d_{k+1}-\psi(\varepsilon),c_{k+1}+\psi(\varepsilon)]$.
\end{lemma}

\section{Proofs}\label{section 3}
Throughout this section, to apply the version of Clark's theorem in \cite{Kajikiya2013}, we work in the  subspace $W=:W_{1}\times W_{2}$ of the Orlicz-Sobolev space $W^{1,\Phi_1}(\mathbb{R}^{N})\times W^{1,\Phi_2}(\mathbb{R}^{N})$ with the norm
\begin{eqnarray}\label{3.1.1}
        \|(u,v)\|
  =    \|u\|_{1,\Phi_1}+ \|v\|_{1,\Phi_2}
  =    \|\nabla u\|_{\Phi_1}+\|u\|_{\Phi_1,V_1}
     +\|\nabla v\|_{\Phi_2}+\|v\|_{\Phi_2,V_2},
\end{eqnarray}
where the subspace $W_{1}$ of $W^{1,\Phi_1}(\mathbb{R}^{N})$ is defined by
\begin{eqnarray}\label{3.1.2}
  W_{1}=\left\{u\in W^{1,\Phi_1}(\mathbb{R}^{N})\Big| \int_{\mathbb{R}^N}V_{1}(x)\Phi_{1}(|u|)dx<\infty\right\}
\end{eqnarray}
with the norm
\begin{eqnarray}\label{3.1.3}
  \|u\|_{1,\Phi_1}=\|\nabla u\|_{\Phi_1}+\|u\|_{\Phi_1,V_{1}},
\end{eqnarray}
where
$$
  \|u\|_{\Phi_1,V_{1}}=\inf\left\{\alpha>0\Big|\int_{\mathbb{R}^N} V_{1}(x)\Phi_1\left(\frac{|u|}{\alpha}\right)dx\le 1\right\}
$$
with $\inf_{x\in\mathbb{R}^N} V_{1}(x)>0$, and
the subspace $W_{2}$ of $W^{1,\Phi_2}(\mathbb{R}^{N})$ is defined by
\begin{eqnarray}\label{3.1.4}
  W_{2}=\left\{v\in W^{1,\Phi_2}(\mathbb{R}^{N})\Big| \int_{\mathbb{R}^N}V_{2}(x)\Phi_{2}(|v|)dx<\infty\right\}
\end{eqnarray}
with the norm
\begin{eqnarray}\label{3.1.5}
  \|v\|_{1,\Phi_2}=\|\nabla v\|_{\Phi_2}+\|v\|_{\Phi_2,V_{2}},
\end{eqnarray}
where
$$
  \|v\|_{\Phi_2,V_{2}}=\inf\left\{\beta>0\Big|\int_{\mathbb{R}^N} V_{2}(x)\Phi_2\left(\frac{|v|}{\beta}\right)dx\le 1\right\}
$$
with $\inf_{x\in\mathbb{R}^N} V_{2}(x)>0$.
It is easy to see that  $(W,\|\cdot\|)$ is a separable and reflexive Banach space.
\par
In what follows, we will give some lemmas to prove our main theorem. To this end, we firstly make some modifies.
For some $\delta>0$, let $\rho\in C^1(\mathbb{R}\times \mathbb{R},[0,1])$ be a cut-off even function defined by
\begin{eqnarray}\label{3.1.6}
 \rho(t,s)= \begin{cases}
           1, \;\;\;  \text { if }\;\;|(t,s)| \le \delta/2, \\
           0, \;\;\;  \text { if }\;\;|(t,s)| \ge \delta
        \end{cases}
\end{eqnarray}
and $t\rho_t(t,s)+s\rho_s(t,s)\le 0 $ for all $(t,s)\in \mathbb{R}^2$.
Next we give some properties for modified functions $\widetilde{F}$ and $\widetilde{G}$ defined by $\widetilde{F}(x,u,v)=\rho(u,v)F(x,u,v)$ and $\widetilde{G}(u,v)=\rho(u,v)G(u,v)$, respectively. By the definitions of $\widetilde{F}(x,u,v)$ and $\widetilde{G}(u,v)$, we have
$$
\widetilde{F}_{u}(x,u,v)=\rho_{u}(u,v)F(x,u,v)+\rho(u,v)F_{u}(x,u,v),\;\;
\widetilde{F}_{v}(x,u,v)=\rho_{v}(u,v)F(x,u,v)+\rho(u,v)F_{v}(x,u,v),
$$
$$
\widetilde{G}_{u}(u,v)=\rho_{u}(u,v)G(u,v)+\rho(u,v)G_{u}(u,v),\;\;
\widetilde{G}_{v}(u,v)=\rho_{v}(u,v)G(u,v)+\rho(u,v)G_{v}(u,v).
$$
\par
\noindent
\begin{lemma}\label{lemma3.1}
 $\widetilde{F}$ and $\widetilde{G}$ are continuous functions and the following conditions hold:
\begin{itemize}
\item[$(F_1')$] there exist two  constants $C_1, C_2>0$ such that
\begin{equation*}
 \left\{
  \begin{array}{l}
 |\widetilde{F}_u(x,u,v)|\leq C_1\left(|u|^{l_{1}r-1}+|v|^{\frac{l_{2}(l_{1}r-1)}{l_{1}}}\right),\\
 |\widetilde{F}_v(x,u,v)|\leq C_2\left(|v|^{l_{2}r-1}+|u|^{\frac{l_{1}(l_{2}r-1)}{l_{2}}}\right)
    \end{array}
 \right.
 \end{equation*}
for all  $(x, u,v)\in  \mathbb{R}^{N}\times \mathbb{R}\times \mathbb{R}$, where
$r\in \left(\max\left\{\frac{1}{l_{1}},\frac{1}{l_{2}}\right\}, \min\left\{\frac{l_{1}}{m_{1}},\frac{l_{2}}{m_{2}}\right\}\right)$;
 \item[$(F_2')$] there exists  a constant $\alpha\in \left[\max\{l_{1},l_{2}\}r,\min\{l_{1},l_{2}\}\right)$ such that
{\small$$
      \widetilde{F}_u(x,u,v)u+\widetilde{F}_v(x,u,v)v-\alpha\widetilde{F}(x,u,v)
\leq  a_{1}(x)|u|^{l_{1}}
    +a_{2}(x)|v|^{l_{2}}
$$}
for all  $(x, u,v)\in  \mathbb{R}^{N}\times \mathbb{R}\times \mathbb{R}$, where $a_{1}(x)=\frac{\min\left\{C_{1,1},\Phi_{1}(1)\right\}\left(l_{1}-\alpha\right)}{l_{1}}V_{1}(x)$
and
$a_{2}(x)=\frac{\min\left\{C_{2,1},\Phi_{2}(1)\right\}\left(l_{2}-\alpha\right)}{l_{2}} V_{2}(x)$;
\item[$(F_3')$]
$\lim_{|(u,v)|\rightarrow 0}\left(\inf_{x \in \mathbb{R}^{N}}\frac{\widetilde{F}(x,u,v)}{|u|^{\alpha}+|v|^{\alpha}}\right)=+\infty$;
\item[$(F_4')$]
 $\widetilde{F}(x,-u,-v)=\widetilde{F}(x,u,v)$
for all  $(x, u,v)\in  \mathbb{R}^{N}\times \mathbb{R}\times \mathbb{R}$;
\item[$(F_5')$]
 $\widetilde{F}(x,u,v)=\widetilde{F}_{u}(x,u,v)=\widetilde{F}_{v}(x,u,v)=\widetilde{G}(u,v)=\widetilde{G}_{u}(u,v)=\widetilde{G}_{v}(u,v)=0$
for all $(u,v)\in  \mathbb{R}\times \mathbb{R}$ with $|(u,v)|\geq \delta$  and $x\in  \mathbb{R}^{N}$.
\end{itemize}
\end{lemma}
\vskip2mm
Define the energy functional $I_{\varepsilon}$ on $W$ corresponding to system \eqref{eq1} is
 \begin{eqnarray}\label{3.1.7}
&           I_{\varepsilon}(u,v):
=           \int_{\mathbb{R}^{N}}\Phi_1(|\nabla u|)dx
          + \int_{\mathbb{R}^{N}}\Phi_2(|\nabla v|)dx
          + \int_{\mathbb{R}^{N}}V_{1}(x)\Phi_1(| u|)dx
          + \int_{\mathbb{R}^{N}}V_{2}(x)\Phi_2(| v|)dx
          \nonumber\\
&\qquad\qquad\qquad\qquad\qquad\qquad
          - \int_{\mathbb{R}^{N}}\widetilde{F}(x,u,v)dx
          - \varepsilon\int_{\mathbb{R}^{N}}\widetilde{G}(u,v)dx,\quad (u,v)\in W.
 \end{eqnarray}
 Under the assumptions $(\phi_1)$--$(\phi_4)$, $(F_0)$  and $(F_1')$, by using standard arguments as in \cite{wang2017,Huidobro1999}, we can prove that $I_{\varepsilon}$ is well-defined. Moreover, it is easy to say that $I_{\varepsilon}$ is of class $C^{1}$ functional on $W$.
 \par
\noindent
\begin{remark}\label{remark3.2}
By Remark \ref{remark 2.9}, under conditions $(\phi_1)$--$(\phi_4)$ and $(V)$,
 the following embeddings are continuous:
 $$
     W_{i}\hookrightarrow L^{\Phi_i}(\mathbb{R}^{N})
      ~\mbox{ and }~
     W_{i}\hookrightarrow L^{p_i}(\mathbb{R}^{N}),
     ~~i=1,2,
 $$
  where $p_i \in [l_{i},l_{i}^{\ast}]$
and the embedding
\begin{eqnarray}\label{3.1.8}
W_{i}(B_R)\hookrightarrow L^{p_i}(B_R), ~i=1,2
\end{eqnarray}
is compact, where $R>0$ and $B_R:=\left\{x\in \mathbb{R}^N: |x|<R\right\}$.
 As a result, there exist positive constants $C_{i,3},\;i=1,2,$ such that
\begin{eqnarray}\label{3.1.9}
       \|u\|_{L^{p_{i}}}\leq C_{i,3}\|u\|_{1,\Phi_{i}},
\end{eqnarray}
where $p_i \in [l_{i},l_{i}^{\ast}],\;i=1,2$.
\end{remark}
\par
\noindent
\begin{remark}\label{remark3.3}
 By the Young's inequality, $(F_0)$, $(F_1)'$, $\widetilde{F}(x,0,0)=0$ and the fact
 $$\widetilde{F}(x,u,v)=\int_{0}^{u}\widetilde{F}_s(x,s,v)ds+\int_{0}^{v}\widetilde{F}_t(x,0,t)dt+\widetilde{F}(x,0,0), \quad \forall  (x,u,v)\in\mathbb{R}^{N}\times\mathbb{R}\times\mathbb{R},$$
it is easy to see that there exists a constant $C_4>0$ such that
\begin{equation}\label{3.1.10}
|\widetilde{F}(x,u,v)|\leq C_4(|u|^{l_{1}r}+|v|^{l_{2}r}), \quad \forall  (x,u,v)\in\mathbb{R}^{N}\times\mathbb{R}\times\mathbb{R}.
\end{equation}
\end{remark}
\par
Next, we prove the {\rm(P-S)} conditions for functional $I_{\varepsilon}$. Until then,  we present the following lemmas which will be used in Lemma \ref{lemma3.6}.

\par
\noindent
\begin{lemma}\label{lemma3.4}
Assume that $(\phi_4)$, $(F_1')$ and $(V)$ hold. If there exists a sequence  $(u_{n},v_{n}) \rightharpoonup (u,v)$ in $W$, then
\begin{eqnarray}\label{3.1.11}
    \int_{\mathbb{R}^N}\left|\widetilde{F}_{u}(x,u_n,v_n)-\widetilde{F}_{u}(x,u,v)\right||u_n-u|dx\rightarrow 0
\;\mbox{and}\;
    \int_{\mathbb{R}^N}\left|\widetilde{F}_{v}(x,u_n,v_n)-\widetilde{F}_{v}(x,u,v)\right||v_n-v|dx\rightarrow 0
\end{eqnarray}
as $n\rightarrow \infty$.
\end{lemma}
\par
\noindent
{\bf Proof.} Evidently, it suffices to show that
\begin{eqnarray}\label{3.1.12}
   \int_{\mathbb{R}^{N}\backslash B_{R}(0)}\left|\widetilde{F}_{u}(x,u_n,v_n)-\widetilde{F}_{u}(x,u,v)\right||u_n-u|dx
+  \int_{B_{R}(0)}\left|\widetilde{F}_{u}(x,u_n,v_n)-\widetilde{F}_{u}(x,u,v)\right||u_n-u|dx
=  o(1).
\end{eqnarray}
First,
by condition $(V)$, for any given constant $\epsilon>0$, there exists $R=R(\epsilon)>0$ such that
\begin{eqnarray}\label{3.1.13}
    \left(\int_{\mathbb{R}^{N}\backslash B_{R}(0)}\left(V_{i}(x)^{-1}\right)^{\frac{r}{1-r}}dx\right)^{1-r}
<   \epsilon,\;\;\forall\;i=1,2.
\end{eqnarray}
Moreover, by $(\phi_4)$ and Lemma \ref{lemma2.6}, we have
\begin{eqnarray}\label{3.1.14}
\Phi_{i}(t)\geq C_{i,4}|t|^{l_{i}}, \forall \; t\geq 0,
\end{eqnarray}
where $C_{i,4}=\min\left\{C_{i,1},\Phi_{i}(1)\right\}$.
Let $\mathbb{A}:=\mathbb{R}^{N}\backslash B_{R}(0)$, by Young's inequality and H\"{o}lder's inequality, we get that
\begin{align}\label{0.1.2}
&
\left(
\int_{\mathbb{A}}\left(V_{1}(x)^{-1}\right)^{\frac{1}{l_1(1-r)}}\left(V_{2}(x)^{-1}\right)^{\frac{l_1r-1}{l_1(1-r)}}
 dx
\right)^{1-r}
\leq
\left(
 \int_{\mathbb{A}}
 \left(
  \frac{1}{l_1r}\left(V_{1}(x)^{-1}\right)^{\frac{r}{1-r}}
+ \frac{l_1r-1}{l_1r}\left(V_{2}(x)^{-1}\right)^{\frac{r}{1-r}}
 \right)
 dx
 \right)^{1-r}
      \nonumber \\
&\qquad
\leq
\left(\frac{1}{l_1r}\right)^{1-r}
\left(\int_{\mathbb{A}}\left(V_{1}(x)^{-1}\right)^{\frac{r}{1-r}}dx\right)^{1-r}
+
\left(\frac{l_1r-1}{l_1r}\right)^{1-r}
\left(\int_{\mathbb{A}}\left(V_{2}(x)^{-1}\right)^{\frac{r}{1-r}}dx\right)^{1-r}
<   \epsilon.
 \end{align}
Thus, by $(F_1')$, (\ref{3.1.13}), (\ref{0.1.2}), Young's inequality and H\"{o}lder's inequality together with (\ref{3.1.14}), we get
\begin{align}\label{3.1.15}
&~~~~\int_{\mathbb{A}}\left|\widetilde{F}_{u}(x,u_n,v_n)-\widetilde{F}_{u}(x,u,v)\right||u_n-u|dx
        \nonumber\\
&
\leq
        C_{1}\int_{\mathbb{A}}|u_n|^{l_{1}r}dx
      + C_{1}\int_{\mathbb{A}}|u_n|^{l_{1}r-1}|u|dx
      + C_{1}\int_{\mathbb{A}}|v_n|^{\frac{l_{2}(l_{1}r-1)}{l_{1}}}|u_n|dx
      + C_{1}\int_{\mathbb{A}}|v_n|^{\frac{l_{2}(l_{1}r-1)}{l_{1}}}|u|dx
      \nonumber \\
&\qquad
      + C_{1}\int_{\mathbb{A}}|u|^{l_{1}r-1}|u_{n}|dx
      + C_{1}\int_{\mathbb{A}}|u|^{l_{1}r}dx
      + C_{1}\int_{\mathbb{A}}|v|^{\frac{l_{2}(l_{1}r-1)}{l_{1}}}|u_n|dx
      + C_{1}\int_{\mathbb{A}}|v|^{\frac{l_{2}(l_{1}r-1)}{l_{1}}}|u|dx
        \nonumber\\
&
\leq
C_{1}\left(\int_{\mathbb{A}}\left(V_{1}(x)^{-1}\right)^{\frac{r}{1-r}}dx\right)^{1-r}
\left[
\left(\int_{\mathbb{A}}V_{1}(x)|u_n|^{l_{1}}dx\right)^{r}
+
\left(\int_{\mathbb{A}}V_{1}(x)|u_n|^{l_{1}}dx\right)^{\frac{l_{1}r-1}{l_{1}}}
\left(\int_{\mathbb{A}}V_{1}(x)|u|^{l_{1}}dx\right)^{\frac{1}{l_{1}}}
\right]
      \nonumber \\
& \qquad
+ C_{1}
\left(\int_{\mathbb{A}}
\left[\left(V_{1}(x)^{-1}\right)^{\frac{1}{l_{1}}}\left(V_{2}(x)^{-1}\right)^{\frac{l_{1}r-1}{l_{1}}}\right]^{\frac{1}{1-r}}
                 dx
\right)^{1-r}
\left(\int_{\mathbb{A}}V_{2}(x)|v_n|^{l_{2}}dx\right)^{\frac{l_{1}r-1}{l_{1}}}
\left(\int_{\mathbb{A}}V_{1}(x)|u_n|^{l_{1}}dx\right)^{\frac{1}{l_{1}}}
      \nonumber \\
& \qquad
+ C_{1}
\left(\int_{\mathbb{A}}
\left[\left(V_{1}(x)^{-1}\right)^{\frac{1}{l_{1}}}\left(V_{2}(x)^{-1}\right)^{\frac{l_{1}r-1}{l_{1}}}\right]^{\frac{1}{1-r}}
                 dx
\right)^{1-r}
\left(\int_{\mathbb{A}}V_{2}(x)|v_n|^{l_{2}}dx\right)^{\frac{l_{1}r-1}{l_{1}}}
\left(\int_{\mathbb{A}}V_{1}(x)|u|^{l_{1}}dx\right)^{\frac{1}{l_{1}}}
      \nonumber \\
& \qquad
+ C_{1}
\left(\int_{\mathbb{A}}\left(V_{1}(x)^{-1}\right)^{\frac{r}{1-r}}dx\right)^{1-r}
\left[
\left(\int_{\mathbb{A}}V_{1}(x)|u|^{l_{1}}dx\right)^{\frac{l_{1}r-1}{l_{1}}}
\left(\int_{\mathbb{A}}V_{1}(x)|u_n|^{l_{1}}dx\right)^{\frac{1}{l_{1}}}
+
\left(\int_{\mathbb{A}}V_{1}(x)|u_n|^{l_{1}}dx\right)^{r}
\right]
      \nonumber \\
& \qquad
+ C_{1}\left(\int_{\mathbb{A}}
\left[\left(V_{1}(x)^{-1}\right)^{\frac{1}{l_{1}}}\left(V_{2}(x)^{-1}\right)^{\frac{l_{1}r-1}{l_{1}}}\right]^{\frac{1}{1-r}}
                 dx
\right)^{1-r}
\left(\int_{\mathbb{A}}V_{2}(x)|v|^{l_{2}}dx\right)^{\frac{l_{1}r-1}{l_{1}}}
\left(\int_{\mathbb{A}}V_{1}(x)|u_n|^{l_{1}}dx\right)^{\frac{1}{l_{1}}}
      \nonumber \\
& \qquad
+ C_{1}
\left(\int_{\mathbb{A}}
\left[\left(V_{1}(x)^{-1}\right)^{\frac{1}{l_{1}}}\left(V_{2}(x)^{-1}\right)^{\frac{l_{1}r-1}{l_{1}}}\right]^{\frac{1}{1-r}}
                 dx
\right)^{1-r}
\left(\int_{\mathbb{A}}V_{2}(x)|v|^{l_{2}}dx\right)^{\frac{l_{1}r-1}{l_{1}}}
\left(\int_{\mathbb{A}}V_{1}(x)|u|^{l_{1}}dx\right)^{\frac{1}{l_{1}}}
        \nonumber\\
&
\leq
     \epsilon \frac{C_{1}}{C_{1,4}^{r}}
     \left(\int_{\mathbb{A}}V_{1}(x)\Phi_{1}(|u_n|)dx\right)^{r}
+   \epsilon \frac{C_{1}}{C_{1,4}^{r}}
\left(\int_{\mathbb{A}}V_{1}(x)\Phi_{1}(|u_n|)dx\right)^{\frac{l_{1}r-1}{l_{1}}}
\left(\int_{\mathbb{A}}V_{1}(x)\Phi_{1}(|u|)dx\right)^{\frac{1}{l_{1}}}
      \nonumber \\
& \qquad
+  \frac{C_{1}}{C_{1,4}^{\frac{1}{l_{1}}}C_{2,4}^{\frac{l_{1}r-1}{l_{1}}}}
\left(\int_{\mathbb{A}}
\left[\left(V_{1}(x)^{-1}\right)^{\frac{1}{l_{1}}}\left(V_{2}(x)^{-1}\right)^{\frac{l_{1}r-1}{l_{1}}}\right]^{\frac{1}{1-r}}
                 dx
\right)^{1-r}
\left(\int_{\mathbb{A}}V_{2}(x)\Phi_{2}(|v_n|)dx\right)^{\frac{l_{1}r-1}{l_{1}}}
\left(\int_{\mathbb{A}}V_{1}(x)\Phi_{1}(|u_n|)dx\right)^{\frac{1}{l_{1}}}
      \nonumber \\
& \qquad
+ \frac{C_{1}}{C_{1,4}^{\frac{1}{l_{1}}}C_{2,4}^{\frac{l_{1}r-1}{l_{1}}}}
\left(\int_{\mathbb{A}}
\left[\left(V_{1}(x)^{-1}\right)^{\frac{1}{l_{1}}}\left(V_{2}(x)^{-1}\right)^{\frac{l_{1}r-1}{l_{1}}}\right]^{\frac{1}{1-r}}
                 dx
\right)^{1-r}
\left(\int_{\mathbb{A}}V_{2}(x)\Phi_{2}(|v_n|)dx\right)^{\frac{l_{1}r-1}{l_{1}}}
\left(\int_{\mathbb{A}}V_{1}(x)\Phi_{1}(|u|)dx\right)^{\frac{1}{l_{1}}}
      \nonumber \\
& \qquad
+ \epsilon \frac{C_{1}}{C_{1,4}^{r}}
\left(\int_{\mathbb{A}}V_{1}(x)\Phi_{1}(|u|)dx\right)^{\frac{l_{1}r-1}{l_{1}}}
\left(\int_{\mathbb{A}}V_{1}(x)\Phi_{1}(|u_n|)dx\right)^{\frac{1}{l_{1}}}
+ \epsilon \frac{C_{1}}{C_{1,4}^{r}}
\left(\int_{\mathbb{A}}V_{1}(x)\Phi_{1}(|u|)dx\right)^{r}
      \nonumber \\
& \qquad
+ \frac{C_{1}}{C_{1,4}^{\frac{1}{l_{1}}}C_{2,4}^{\frac{l_{1}r-1}{l_{1}}}}
\left(\int_{\mathbb{A}}
\left[\left(V_{1}(x)^{-1}\right)^{\frac{1}{l_{1}}}\left(V_{2}(x)^{-1}\right)^{\frac{l_{1}r-1}{l_{1}}}\right]^{\frac{1}{1-r}}
                 dx
\right)^{1-r}
\left(\int_{\mathbb{A}}V_{2}(x)\Phi_{2}(|v|)dx\right)^{\frac{l_{1}r-1}{l_{1}}}
\left(\int_{\mathbb{A}}V_{1}(x)\Phi_{1}(|u_n|)dx\right)^{\frac{1}{l_{1}}}
      \nonumber \\
& \qquad
+ \frac{C_{1}}{C_{1,4}^{\frac{1}{l_{1}}}C_{2,4}^{\frac{l_{1}r-1}{l_{1}}}}
\left(\int_{\mathbb{A}}
\left[\left(V_{1}(x)^{-1}\right)^{\frac{1}{l_{1}}}\left(V_{2}(x)^{-1}\right)^{\frac{l_{1}r-1}{l_{1}}}\right]^{\frac{1}{1-r}}
                 dx
\right)^{1-r}
\left(\int_{\mathbb{A}}V_{2}(x)\Phi_{2}(|v|)dx\right)^{\frac{l_{1}r-1}{l_{1}}}
\left(\int_{\mathbb{A}}V_{1}(x)\Phi_{1}(|u|)dx\right)^{\frac{1}{l_{1}}}
        \nonumber\\
&
<
\epsilon \frac{C_{1}}{C_{1,4}^{r}}
\left(
  \max\left\{\|u_n\|_{1,\Phi_{1}}^{l_{1}r},\|u_n\|_{1,\Phi_{1}}^{m_{1}r}\right\}
+
\max\left\{\|u_n\|_{1,\Phi_{1}}^{l_{1}r-1},\|u_n\|_{1,\Phi_{1}}^{\frac{m_{1}(l_{1}r-1)}{l_{1}}}\right\}
\max\left\{\|u\|_{1,\Phi_{1}},\|u\|_{1,\Phi_{1}}^{\frac{m_{1}}{l_{1}}}\right\}
\right)
      \nonumber \\
& \qquad
      +  \epsilon\frac{C_{1}}{C_{1,4}^{\frac{1}{l_{1}}}C_{2,4}^{\frac{l_{1}r-1}{l_{1}}}}
        \max\left\{\|v_n\|_{1,\Phi_{2}}^{\frac{l_{2}(l_{1}r-1)}{l_{1}}},\|v_n\|_{1,\Phi_{2}}^{\frac{m_{2}(l_{1}r-1)}{l_{1}}}\right\}
        \max\left\{\|u_{n}\|_{1,\Phi_{1}},\|u_{n}\|_{1,\Phi_{1}}^{\frac{m_{1}}{l_{1}}}\right\}
      \nonumber \\
& \qquad
      + \epsilon\frac{C_{1}}{C_{1,4}^{\frac{1}{l_{1}}}C_{2,4}^{\frac{l_{1}r-1}{l_{1}}}}
        \max\left\{\|v_n\|_{1,\Phi_{2}}^{\frac{l_{2}(l_{1}r-1)}{l_{1}}},\|v_n\|_{1,\Phi_{2}}^{\frac{m_{2}(l_{1}r-1)}{l_{1}}}\right\}
        \max\left\{\|u\|_{1,\Phi_{1}},\|u\|_{1,\Phi_{1}}^{\frac{m_{1}}{l_{1}}}\right\}
      \nonumber \\
& \qquad
      + \epsilon \frac{C_{1}}{C_{1,4}^{r}}
      \left(
        \max\left\{\|u\|_{1,\Phi_{1}}^{l_{1}r-1},\|u\|_{1,\Phi_{1}}^{\frac{m_{1}(l_{1}r-1)}{l_{1}}}\right\}
        \max\left\{\|u_n\|_{1,\Phi_{1}},\|u_n\|_{1,\Phi_{1}}^{\frac{m_{1}}{l_{1}}}\right\}
      +
         \max\left\{\|u\|_{1,\Phi_{1}}^{l_{1}r},\|u\|_{1,\Phi_{1}}^{m_{1}r}\right\}
        \right)
      \nonumber \\
& \qquad
      + \epsilon\frac{C_{1}}{C_{1,4}^{\frac{1}{l_{1}}}C_{2,4}^{\frac{l_{1}r-1}{l_{1}}}}
        \max\left\{\|v\|_{1,\Phi_{2}}^{\frac{l_{2}(l_{1}r-1)}{l_{1}}},\|v\|_{1,\Phi_{2}}^{\frac{m_{2}(l_{1}r-1)}{l_{1}}}\right\}
        \max\left\{\|u_{n}\|_{1,\Phi_{1}},\|u_{n}\|_{1,\Phi_{1}}^{\frac{m_{1}}{l_{1}}}\right\}
      \nonumber \\
& \qquad
      + \epsilon\frac{C_{1}}{C_{1,4}^{\frac{1}{l_{1}}}C_{2,4}^{\frac{l_{1}r-1}{l_{1}}}}
        \max\left\{\|v\|_{1,\Phi_{2}}^{\frac{l_{2}(l_{1}r-1)}{l_{1}}},\|v\|_{1,\Phi_{2}}^{\frac{m_{2}(l_{1}r-1)}{l_{1}}}\right\}
        \max\left\{\|u\|_{1,\Phi_{1}},\|u\|_{1,\Phi_{1}}^{\frac{m_{1}}{l_{1}}}\right\}.
\end{align}
As $\epsilon$ is arbitrary, \eqref{3.1.15} implies that the first term in equality \eqref{3.1.12} converges to $0$.
Next, we claim that
\begin{eqnarray}\label{3.2.1}
\int_{B_{R}(0)}\left|\widetilde{F}_{u}(x,u_n,v_n)-\widetilde{F}_{u}(x,u,v)\right||u_n-u|dx
=  o(1).
\end{eqnarray}
In fact by $u_{u}\rightharpoonup u$ in $W^{1,\Phi_{1}}(\mathbb{R}^{N})$, it is easy to see that  $u_n|_{B_{R}}\rightarrow u|_{B_{R}}$ in $L^{l_1}(B_{R})$ from Remark \ref{remark3.2}. So for $\epsilon>0$ above, there exists $n_{0}\in \mathbb{N}$ such that
\begin{eqnarray}\label{3.1.17}
\left(\int_{B_{R}(0)}|u_n-u|^{l_{1}}dx\right)^{\frac{1}{l_{1}}}<\epsilon,\;\;\forall\; n>n_{0}.
\end{eqnarray}
Then for all $n>n_{0}$, by $(F_1')$, Young's inequality and H\"{o}lder's inequality,
 we have
\begin{eqnarray}\label{3.1.18}
&&\int_{B_{R}(0)}\left|\widetilde{F}_{u}(x,u_n,v_n)-\widetilde{F}_{u}(x,u,v)\right||u_n-u|dx
        \nonumber\\
&&
 \leq   C_{1}\left(\int_{B_{R}(0)}
                   \left(|u_n|^{l_{1}r-1}+|u|^{l_{1}r-1}+|v_n|^{\frac{l_{2}(l_{1}r-1)}{l_{1}}}+|v|^{\frac{l_{2}(l_{1}r-1)}{l_{1}}}\right)^{\frac{l_{1}}{l_{1}-1}}
                   dx
           \right)^{\frac{l_{1}-1}{l_{1}}}
           \left(\int_{B_{R}(0)}|u_n-u|^{l_{1}}dx\right)^{\frac{1}{l_{1}}}
        \nonumber\\
&&
 \leq   \epsilon C_{1}
        \left(\int_{B_{R}(0)}
                \left(
                   2^{\frac{1}{l_{1}-1}}
                   \left(   |u_n|^{l_{1}r-1}
                          + |u|^{l_{1}r-1}
                   \right)^{\frac{l_{1}}{l_{1}-1}}
         +         2^{\frac{1}{l_{1}-1}}
                   \left(   |v_n|^{\frac{l_{2}(l_{1}r-1)}{l_{1}}}
                          + |v|^{\frac{l_{2}(l_{1}r-1)}{l_{1}}}
                   \right)^{\frac{l_{1}}{l_{1}-1}}
               \right)
                   dx
           \right)^{\frac{l_{1}-1}{l_{1}}}
        \nonumber\\
&&
 \leq   \epsilon C_{1} 2^{\frac{1}{l_{1}}}
        \left(\int_{B_{R}(0)}
                \left(
                   |u_n|^{l_{1}}
                 + |u|^{l_{1}}
                 + |v_n|^{l_{2}}
                 + |v|^{l_{2}}
                 + 4
              \right)
                   dx
           \right)^{\frac{l_{1}-1}{l_{1}}}
        \nonumber\\
&&
=      \epsilon C_{1} 2^{\frac{1}{l_{1}}}
            \left(
                   \|u_n\|_{L^{l_{1}}}^{l_{1}}
                 + \|u\|_{L^{l_{1}}}^{l_{1}}
                 + \|v_n\|_{L^{l_{2}}}^{l_{2}}
                 + \|v\|_{L^{l_{2}}}^{l_{2}}
                 + 4\mu(B_{R}(0))
           \right)^{\frac{l_{1}-1}{l_{1}}},
\end{eqnarray}
where $\mu(B_{R}(0))$ denotes the Lebesgue measure of $B_{R}(0)$.
Hence, the arbitrariness of $\epsilon$, the boundedness of $\|u_n\|_{L^{l_{1}}}$ and $\|v_n\|_{L^{l_{2}}}$ and \eqref{3.1.18} imply that \eqref{3.2.1} hold. Thus, we have
\begin{eqnarray*}
    \int_{\mathbb{R}^N}\left|\widetilde{F}_{u}(x,u_n,v_n)-\widetilde{F}_{u}(x,u,v)\right||u_n-u|dx=o(1).
\end{eqnarray*}
Similarly, we can obtain that
\begin{eqnarray*}
    \int_{\mathbb{R}^N}\left|\widetilde{F}_{v}(x,u_n,v_n)-\widetilde{F}_{v}(x,u,v)\right||u_n-u|dx=o(1).
\end{eqnarray*}
The proof is complete.\qed
\par
\noindent
\begin{lemma}\label{lemma3.5}
Assume that $(K)$ and $(G)$ hold. If there exists a sequence  $(u_{n},v_{n}) \rightharpoonup (u,v)$ in $W$, then
\begin{eqnarray}\label{3.1.19}
    \int_{\mathbb{R}^N}k(x)\left|\widetilde{G}_{u}(u_n,v_n)-\widetilde{G}_{u}(u,v)\right||u_n-u|dx\rightarrow 0
\;\mbox{and}\;
    \int_{\mathbb{R}^N}k(x)\left|\widetilde{G}_{v}(u_n,v_n)-\widetilde{G}_{v}(u,v)\right||v_n-v|dx\rightarrow 0
\end{eqnarray}
as $n\rightarrow \infty$.
\end{lemma}
\par
\noindent
{\bf Proof.}
Let $h_n, h: \mathbb{R}^N\rightarrow\mathbb{R}$ be defined by
 \begin{eqnarray*}
h_n=k(x)\left|\widetilde{G}_{u}(u_n,v_n)-\widetilde{G}_{u}(u,v)\right||u_n-u|,\;\;\;\;\;\;
h=4|k(x)| \left|\widetilde{G}_{u}(u,v)\right||u|.
 \end{eqnarray*}
 By  $ u_{n}  \rightarrow  u $ a.e. in $x\in \mathbb{R}^N$ and the continuity of $\widetilde{G}_{u}(u,v)$, we get
\begin{eqnarray}\label{3.1.20}
 |h_n| &  =   & |k(x)|\left|\widetilde{G}_{u}(u_n,v_n)-\widetilde{G}_{u}(u,v)\right||u_n-u|\nonumber\\
       & \leq & |k(x)|\left(   \left|\widetilde{G}_{u}(u_n,v_n)\right||u_n|
                                             + \left|\widetilde{G}_{u}(u_n,v_n)\right||u|
                                             + \left|\widetilde{G}_{u}(u,v)\right||u_n|
                                             + \left|\widetilde{G}_{u}(u,v)\right||u|
                                             \right)\nonumber\\
      &  \rightarrow   & 4|k(x)| \left|\widetilde{G}_{u}(u,v)\right||u|=h.
\end{eqnarray}
Recall that $\widetilde{G}_{u}(u,v)=0$ for $|(u,v)|\geq \delta$. So $h\in L^{1}(\mathbb{R}^N)$. From this and Lebesgue dominated convergence theorem, we have
\begin{eqnarray*}
\lim_{n\rightarrow\infty}\int_{\mathbb{R}^N}k(x)\left|\widetilde{G}_{u}(u_n,v_n)-\widetilde{G}_{u}(u,v)\right||u_n-u|dx=0.
\end{eqnarray*}
Similarly, we can obtain that
\begin{eqnarray*}
\lim_{n\rightarrow\infty}\int_{\mathbb{R}^N}k(x)\left|\widetilde{G}_{v}(u_n,v_n)-\widetilde{G}_{v}(u,v)\right||v_n-v|dx=0.
\end{eqnarray*}
The proof is complete.\qed
\par
\noindent
\begin{lemma}\label{lemma3.6}
Assume that $(\phi_1)$-$(\phi_4)$, $(F_1')$ and $(V)$ hold.  For $\varepsilon\in [0,1]$, $I_{\varepsilon}$  satisfies {\rm(P-S)} conditions uniformly on $\varepsilon$.
\end{lemma}
\par
\noindent
{\bf Proof.}
First, we claim that the functional $I_{\varepsilon}(u,v)$ is coercive and bounded below. In fact, for any given $(u,v)\in W$, by $(\phi_4)$, $(V)$, \eqref{3.1.10}  and  H\"{o}lder's inequality, we have
\begin{align}\label{3.1.21}
&\int_{\mathbb{R}^{N}} |\widetilde{F}(x,u,v)|dx
 \nonumber\\
&\quad
=
        C_{4}\int_{\mathbb{R}^{N}} \left(V_{1}^{-1}(x)\right)^{r}
                                   \left(V_{1}(x)\right)^{r}
                                   |u|^{l_{1}r}
             dx
      + C_{4}\int_{\mathbb{R}^{N}} \left(V_{2}^{-1}(x)\right)^{r}
                                   \left(V_{2}(x)\right)^{r}
                                   |v|^{l_{2}r}
              dx\nonumber\\
&\quad
\leq
        C_{4}\left(\int_{\mathbb{R}^{N}} \left(V_{1}^{-1}(x)\right)^{\frac{r}{1-r}}
                    dx
             \right)^{1-r}
             \left(\int_{\mathbb{R}^{N}} V_{1}(x)
                                   |u|^{l_{1}}
                    dx
             \right)^{r}
      + C_{4}\left(\int_{\mathbb{R}^{N}} \left(V_{2}^{-1}(x)\right)^{\frac{r}{1-r}}
                    dx
             \right)^{1-r}
             \left(\int_{\mathbb{R}^{N}} V_{2}(x)
                                   |v|^{l_{2}}
                    dx
             \right)^{r}\nonumber\\
&\quad
\leq
       \frac{C_{4}}{C_{1,4}^{r}}
             \left(\int_{\mathbb{R}^{N}} \left(V_{1}^{-1}(x)\right)^{\frac{r}{1-r}}
                    dx
             \right)^{1-r}
             \left(\int_{\mathbb{R}^{N}} V_{1}(x)
                                         \Phi_1(|u|)
                    dx
             \right)^{r}
      +     \frac{C_{4}}{C_{2,4}^{r}}
             \left(\int_{\mathbb{R}^{N}} \left(V_{2}^{-1}(x)\right)^{\frac{r}{1-r}}
                    dx
             \right)^{1-r}
             \left(\int_{\mathbb{R}^{N}} V_{2}(x)
                                         \Phi_2(|v|)
                    dx
             \right)^{r}
             \nonumber\\
&\quad
\leq
        C_{5}\max\left\{\|u\|_{\Phi_1,V_{1}}^{l_{1}r},\|u\|_{\Phi_1,V_{1}}^{m_{1}r}\right\}
      + C_{5}\max\left\{\|v\|_{\Phi_2,,V_{2}}^{l_{2}r},\|v\|_{\Phi_2,V_{2}}^{m_{2}r}\right\}
             \nonumber\\
&\quad
\leq
        C_{5}\max\left\{\|u\|_{1,\Phi_1}^{l_{1}r},\|u\|_{1,\Phi_1}^{m_{1}r}\right\}
      + C_{5}\max\left\{\|v\|_{1,\Phi_2}^{l_{2}r},\|v\|_{1,\Phi_2}^{m_{2}r}\right\},
\end{align}
where $C_{5}
=\max\left\{\frac{C_{4}}{C_{1,4}^{r}}\|V_{1}^{-1}(x)\|_{L^{\frac{r}{1-r}}}^{r},
\frac{C_{4}}{C_{2,4}^{r}}\|V_{2}^{-1}(x)\|_{L^{\frac{r}{1-r}}}^{r}\right\}
$.
 Moreover, by $(K)$, for $\varepsilon\in [0,1]$, we have
\begin{align}\label{3.1.22}
&    \left|\varepsilon \int_{\mathbb{R}^{N}} k(x)\widetilde{G}(u,v)dx \right|
\leq
        |\varepsilon|\int_{|(u,v)|\leq\delta}|k(x)||\widetilde{G}(u,v)|dx
      + |\varepsilon|\int_{|(u,v)|>\delta} |k(x)||\widetilde{G}(u,v)|dx\nonumber\\
&\qquad\qquad\qquad\quad\quad~~
=
        |\varepsilon|\int_{|(u,v)|\leq\delta} |k(x)||\widetilde{G}(u,v)|dx
\leq
        |\varepsilon|\max_{|(u,v)|\leq\delta}|\widetilde{G}(u,v)|\int_{\mathbb{R}^{N}} |k(x)|dx
\leq
        |\varepsilon|C_{6}.
\end{align}
 Thus,
\begin{align}\label{3.1.23}
&I_{\varepsilon}(u,v)
=         \int_{\mathbb{R}^{N}}\Phi_1(|\nabla u|)dx
          + \int_{\mathbb{R}^{N}}\Phi_2(|\nabla v|)dx
          + \int_{\mathbb{R}^{N}}V_{1}(x)\Phi_1(| u|)dx
          + \int_{\mathbb{R}^{N}}V_{2}(x)\Phi_2(| v|)dx
          \nonumber\\
&~~~~
          - \int_{\mathbb{R}^{N}}\widetilde{F}(x,u,v)dx
          - \varepsilon\int_{\mathbb{R}^{N}}k(x)\widetilde{G}(u,v)dx
          \nonumber\\
&\geq        \min\left\{\|\nabla u\|_{\Phi_1}^{l_{1}},\|\nabla u\|_{\Phi_1}^{m_{1}}\right\}
         + \min\left\{\|u\|_{\Phi_1,V_{1}}^{l_{1}},\|u\|_{\Phi_1,V_{1}}^{m_{1}}\right\}
          - C_{5}\max\left\{\|u\|_{1,\Phi_1}^{l_{1}r},\|u\|_{1,\Phi_1}^{m_{1}r}\right\}
          \nonumber\\
&~~~~
          + \min\left\{\|\nabla v\|_{\Phi_2}^{l_{2}},\|\nabla v\|_{\Phi_2}^{m_{2}}\right\}
         + \min\left\{\|v\|_{\Phi_2,V_{2}}^{l_{2}},\|v\|_{\Phi_2,V_{2}}^{m_{2}}\right\}
          - C_{5}\max\left\{\|v\|_{1,\Phi_2}^{l_{2}r},\|v\|_{1,\Phi_2}^{m_{2}r}\right\}
          - |\varepsilon|C_{6}
          \nonumber\\
&\geq       \|\nabla u\|_{\Phi_1}^{l_{1}}
         + \|u\|_{\Phi_1,V_{1}}^{l_{1}}
         -2
          - C_{5}\max\left\{\|u\|_{1,\Phi_1}^{l_{1}r},\|u\|_{1,\Phi_1}^{m_{1}r}\right\}
          \nonumber\\
&~~~~
          + \|\nabla v\|_{\Phi_2}^{l_{2}}
         +\|v\|_{\Phi_2,V_{2}}^{l_{2}}
         -2
          - C_{5}\max\left\{\|v\|_{1,\Phi_2}^{l_{2}r},\|v\|_{1,\Phi_2}^{m_{2}r}\right\}
          - |\varepsilon|C_{6}
          \nonumber\\
&\geq        \frac{1}{2^{l_{1}-1}}\|u\|_{1,\Phi_1}^{l_{1}}
          - C_{5}\max\left\{\|u\|_{1,\Phi_1}^{l_{1}r},\|u\|_{1,\Phi_1}^{m_{1}r}\right\}
          + \frac{1}{2^{l_{2}-1}}\|v\|_{1,\Phi_2}^{l_{2}}
          - C_{5}\max\left\{\|v\|_{1,\Phi_2}^{l_{2}r},\|v\|_{1,\Phi_2}^{m_{2}r}\right\}
          - |\varepsilon|C_{6}-4.
\end{align}
Now, let $(\varepsilon_{n},(u_{n},v_{n}))\in [0,1]\times W$ be any sequence such that
\begin{equation}\label{3.1.24}
I_{\varepsilon_{n}}(u_{n},v_{n})<\infty \;\mbox{and}\; I_{\varepsilon_{n}}^{'}(u_{n},v_{n})\rightarrow 0.
\end{equation}
Then by \eqref{3.1.23}, $\{\varepsilon_{n}\}$ and $\{(u_{n},v_{n})\}$ are bounded.
Therefore, there exists a subsequence of $\{\varepsilon_{n}\}$, still denoted by $\{\varepsilon_{n}\}$, such that $\{\varepsilon_{n}\}$ converges to $\varepsilon$, and a subsequence of $\{(u_{n},v_{n})\}$, still denoted by $\{(u_{n},v_{n})\}$,  such that $\{(u_{n},v_{n})\}$ converges to $(u,v)$
weakly in $W$ and a.e. on $\R^{N}$.
Next, we shall show that this convergence is a strong convergence actually.
Note that
\begin{eqnarray*}\label{3.1.25}
&&          o(1)
 =    \langle I_{\varepsilon}'(u_{n},v_{n})-I_{\varepsilon}'(u,v),(u_n-u,v_n-v)\rangle
        \nonumber\\
&&\qquad =        \int_{\mathbb{R}^{N}}\left(\phi_1(|\nabla u_{n}|)\nabla u_{n}-\phi_1(|\nabla u|)\nabla u \right)\left(\nabla u_{n}-\nabla u\right)dx
          + \int_{\mathbb{R}^{N}}\left(\phi_2(|\nabla v_{n}|)\nabla v_{n}-\phi_2(|\nabla v|)\nabla v \right)\left(\nabla v_{n}-\nabla v\right)dx
          \nonumber\\
&&\qquad\quad
          + \int_{\mathbb{R}^{N}}V_{1}(x)\left(\phi_1(|u_{n}|)u_{n}-\phi_1(|u|)u \right)\left(u_{n}-u\right)dx
          + \int_{\mathbb{R}^{N}}V_{2}(x)\left(\phi_2(|v_{n}|)v_{n}-\phi_2(|v|)v \right)\left(v_{n}-v\right)dx
          \nonumber\\
&&\qquad\quad
          - \int_{\mathbb{R}^{N}}\left(\widetilde{F}_{u_{n}}(x,u_{n},v_{n})-\widetilde{F}_{u}(x,u,v)\right)(u_{n}-u)dx
          - \int_{\mathbb{R}^{N}}\left(\widetilde{F}_{v_{n}}(x,u_{n},v_{n})-\widetilde{F}_{v}(x,u,v)\right)(v_{n}-v)dx
          \nonumber\\
&&\qquad\quad
 - \varepsilon\int_{\mathbb{R}^{N}}k(x)\left(\widetilde{G}_{u_{n}}(u_{n},v_{n})-\widetilde{G}_{u}(u,v)\right)(u_{n}-u)dx
 - \varepsilon\int_{\mathbb{R}^{N}}k(x)\left(\widetilde{G}_{v_{n}}(u_{n},v_{n})-\widetilde{G}_{v}(u,v)\right)(v_{n}-v)dx.
\end{eqnarray*}
Then by Lemma 3.9 in \cite{Wang2016}, Lemma 3.4 and Lemma 3.5 imply that $\|u_{n}-u\|_{1,\Phi_1}\rightarrow 0$ and $\|v_{n}-v\|_{1,\Phi_2}\rightarrow 0$ as $n\rightarrow\infty$. \qed
\par
\noindent
\begin{lemma}\label{lemma3.7}
Assume that $(\phi_1)$-$(\phi_4)$, $(F_2')$ and $(F_3')$ hold.  For any $(u,v)\in W \backslash \{(0,0)\}$, $I_{0}(u,v)$  satisfies condition $(A_{5})$.
\end{lemma}
\par
\noindent
{\bf Proof.}
For any given $(u,v)\in W \backslash \{(0,0)\}$.
Note that $\widetilde{F}$ is even on $(u,v)\in \mathbb{R} \times \mathbb{R}$. Without loss of generality, we may assume $t>0$ and consider
 \begin{eqnarray}\label{3.1.26}
&&\gamma_{(u,v)}(t)=t^{-\alpha}I_{0}(tu,tv)
=         t^{-\alpha}\int_{\mathbb{R}^{N}}\Phi_1(t|\nabla u|)dx
          + t^{-\alpha}\int_{\mathbb{R}^{N}}\Phi_2(t|\nabla v|)dx
          + t^{-\alpha}\int_{\mathbb{R}^{N}}V_{1}(x)\Phi_1(t| u|)dx
          \nonumber\\
&&\qquad\qquad\qquad\qquad\qquad
          + t^{-\alpha}\int_{\mathbb{R}^{N}}V_{2}(x)\Phi_2(t| v|)dx
          - t^{-\alpha}\int_{\mathbb{R}^{N}}\widetilde{F}(x,tu,tv)dx.
\end{eqnarray}
Then, by $(\phi_3)$ and $(F_2')$, we have
\begin{eqnarray}\label{3.1.27}
&&\gamma_{(u,v)}'(t)
=         t^{-\alpha-1}\int_{\mathbb{R}^{N}}\left(\phi_1(t|\nabla u|)t^{2}|\nabla u|^{2}-\alpha\Phi_1(t|\nabla u|)\right)dx
          + t^{-\alpha-1}\int_{\mathbb{R}^{N}}\left(\phi_2(t|\nabla v|)t^{2}|\nabla v|^{2}-\alpha\Phi_2(t|\nabla v|)\right)dx
          \nonumber\\
&&\qquad\qquad
 + t^{-\alpha-1}\int_{\mathbb{R}^{N}}V_{1}(x)\left(\phi_1(t| u|)t^{2}| u|^{2}-\alpha\Phi_1(t| u|)\right)dx
          + t^{-\alpha-1}\int_{\mathbb{R}^{N}}V_{2}(x)\left(\phi_2(t| v|)t^{2}| v|^{2}-\alpha\Phi_2(t| v|)\right)dx
\nonumber\\
&&\qquad\qquad
          - t^{-\alpha-1}\int_{\mathbb{R}^{N}}\widetilde{F}_{u}(x,tu,tv)tudx
          - t^{-\alpha-1}\int_{\mathbb{R}^{N}}\widetilde{F}_{v}(x,tu,tv)tvdx
          + t^{-\alpha-1}\int_{\mathbb{R}^{N}}\alpha\widetilde{F}(x,tu,tv)dx
          \nonumber\\
&&\qquad\quad\geq      t^{-\alpha-1}\int_{\mathbb{R}^{N}}\left(l_{1}-\alpha\right)\Phi_1(t|\nabla u|)dx
          + t^{-\alpha-1}\int_{\mathbb{R}^{N}}\left(l_{2}-\alpha\right)\Phi_2(t|\nabla v|)dx
          \nonumber\\
&&\qquad\qquad
          + t^{-\alpha-1}\int_{\mathbb{R}^{N}}V_{1}(x)\left(l_{1}-\alpha\right)\Phi_1(t|u|)dx
          + t^{-\alpha-1}\int_{\mathbb{R}^{N}}V_{2}(x)\left(l_{2}-\alpha\right)\Phi_2(t| v|)dx
          \nonumber\\
&&\qquad\qquad
          - t^{-\alpha-1}\int_{\mathbb{R}^{N}}\left(   \widetilde{F}_{u}(x,tu,tv)tu
                                             + \widetilde{F}_{v}(x,tu,tv)tv
                                             - \alpha\widetilde{F}(x,tu,tv)
                                      \right)dx
          \nonumber\\
&&\qquad\quad\geq      t^{-\alpha-1}\int_{\mathbb{R}^{N}}\left(l_{1}-\alpha\right)\Phi_1(t|\nabla u|)dx
          + t^{-\alpha-1}\int_{\mathbb{R}^{N}}\left(l_{2}-\alpha\right)\Phi_2(t|\nabla v|)dx
          \nonumber\\
          \nonumber\\
&&\qquad\qquad
          + t^{-\alpha-1}
          \int_{\mathbb{R}^{N}}
          \left[V_{1}(x)\left(l_{1}-\alpha\right)C_{1,4}|tu|^{l_{1}}
          + V_{2}(x)\left(l_{2}-\alpha\right)C_{2,4}|tv|^{l_{2}}\right]dx
          \nonumber\\
&&\qquad\qquad
           - t^{-\alpha-1}
           \int_{\mathbb{R}^{N}}
           \left(
                    V_{1}(x)\left(l_{1}-\alpha\right)\frac{C_{1,4}}{l_{1}}|tu|^{l_{1}}
                 + V_{2}(x)\left(l_{2}-\alpha\right)\frac{C_{2,4}}{l_{2}}|tv|^{l_{2}}
                                      \right)dx
          \nonumber\\
&&\qquad\quad>       0.
\end{eqnarray}
Take $\eta>0$ small enough such that
$$
\mu(K_{\eta})> 0\;\; \mbox{and}\;\; K_{\eta}:=\left\{x \in \mathbb{R}^{N}: \eta <|(u,v)|<\frac{1}{\eta}\right\},
$$
where $\mu$ denotes the Lebesgue measure of $\mathbb{R}^{N}$.
Due to $\widetilde{F}(x,u,v)\geq 0$, by Lemma \ref{lemma2.6}, we can estimate the function $\gamma_{(u,v)}(t)$ as follows:
for $0<t<1$,
\begin{eqnarray}\label{3.1.28}
&&\gamma_{(u,v)}(t)=t^{-\alpha}I_{0}(tu,tv)
=          t^{-\alpha}\int_{\mathbb{R}^{N}}\Phi_1(t|\nabla u|)dx
          + t^{-\alpha}\int_{\mathbb{R}^{N}}\Phi_2(t|\nabla v|)dx
          + t^{-\alpha}\int_{\mathbb{R}^{N}}V_{1}(x)\Phi_1(t| u|)dx
          \nonumber\\
&&\qquad\qquad
          + t^{-\alpha}\int_{\mathbb{R}^{N}}V_{2}(x)\Phi_2(t| v|)dx
          - t^{-\alpha}\int_{\mathbb{R}^{N}}\widetilde{F}(x,tu,tv)dx
          \nonumber\\
&&\qquad\quad
\leq       t^{l_{1}-\alpha}\int_{\mathbb{R}^{N}}\Phi_1(|\nabla u|)dx
          + t^{l_{2}-\alpha}\int_{\mathbb{R}^{N}}\Phi_2(|\nabla v|)dx
          + t^{l_{1}-\alpha}\int_{\mathbb{R}^{N}}V_{1}(x)\Phi_1(| u|)dx
          \nonumber\\
&&\qquad\qquad
          + t^{l_{2}-\alpha}\int_{\mathbb{R}^{N}}V_{2}(x)\Phi_2(| v|)dx
          - t^{-\alpha}\int_{K_{\eta}}\widetilde{F}(x,tu,tv)dx
          \nonumber\\
&&\qquad\quad
\leq       t^{l_{1}-\alpha}\int_{\mathbb{R}^{N}}\Phi_1(|\nabla u|)dx
          + t^{l_{2}-\alpha}\int_{\mathbb{R}^{N}}\Phi_2(|\nabla v|)dx
          + t^{l_{1}-\alpha}\int_{\mathbb{R}^{N}}V_{1}(x)\Phi_1(| u|)dx
          \nonumber\\
&&\qquad\qquad
          + t^{l_{2}-\alpha}\int_{\mathbb{R}^{N}}V_{2}(x)\Phi_2(| v|)dx
          - \frac{1}{C_{\alpha}}\int_{K_{\eta}}\frac{\widetilde{F}(x,tu,tv)(|u|^{2}+|v|^{2})^{\frac{\alpha}{2}}}{|tu|^{\alpha}+|tv|^{\alpha}}dx
          \nonumber\\
&&\qquad\quad
\leq       t^{l_{1}-\alpha}\int_{\mathbb{R}^{N}}\Phi_1(|\nabla u|)dx
          + t^{l_{2}-\alpha}\int_{\mathbb{R}^{N}}\Phi_2(|\nabla v|)dx
          + t^{l_{1}-\alpha}\int_{\mathbb{R}^{N}}V_{1}(x)\Phi_1(| u|)dx
          \nonumber\\
&&\qquad\qquad
          + t^{l_{2}-\alpha}\int_{\mathbb{R}^{N}}V_{2}(x)\Phi_2(| v|)dx
          - \frac{\eta^{\alpha}}{C_{\alpha}}\int_{K_{\eta}}\frac{\widetilde{F}(x,tu,tv)}{|tu|^{\alpha}+|tv|^{\alpha}}dx
          \nonumber\\
&&\qquad\quad
\leq       t^{l_{1}-\alpha}\int_{\mathbb{R}^{N}}\Phi_1(|\nabla u|)dx
          + t^{l_{2}-\alpha}\int_{\mathbb{R}^{N}}\Phi_2(|\nabla v|)dx
          + t^{l_{1}-\alpha}\int_{\mathbb{R}^{N}}V_{1}(x)\Phi_1(| u|)dx
          \nonumber\\
&&\qquad\qquad
          + t^{l_{2}-\alpha}\int_{\mathbb{R}^{N}}V_{2}(x)\Phi_2(| v|)dx
          - \frac{\eta^{\alpha}}{C_{\alpha}}\mu(K_{\eta})\inf_{x\in K_{\eta}}\frac{\widetilde{F}(x,tu,tv)}{|tu|^{\alpha}+|tv|^{\alpha}}.
\end{eqnarray}
By condition $(F_3')$, it implies that $\lim_{t\rightarrow 0^{+}}\gamma_{(u,v)}(t)=-\infty$. Hence,  $\gamma_{(u,v)}(t)<0$ for $t>0$ small enough.
\par
On the other hand, by \eqref{3.1.21} and Lemma \ref{lemma2.6}, we have, for $t>0$ large enough,
\begin{eqnarray}\label{3.1.29}
&&\gamma_{(u,v)}(t)=t^{-\alpha}I_{0}(tu,tv)
=          t^{-\alpha}\int_{\mathbb{R}^{N}}\Phi_1(t|\nabla u|)dx
          + t^{-\alpha}\int_{\mathbb{R}^{N}}\Phi_2(t|\nabla v|)dx
          + t^{-\alpha}\int_{\mathbb{R}^{N}}V_{1}(x)\Phi_1(t| u|)dx
          \nonumber\\
&&\qquad\qquad
          + t^{-\alpha}\int_{\mathbb{R}^{N}}V_{2}(x)\Phi_2(t| v|)dx
          - t^{-\alpha}\int_{\mathbb{R}^{N}}\widetilde{F}(x,tu,tv)dx
          \nonumber\\
&&\qquad\quad
\geq       t^{l_{1}-\alpha}\int_{\mathbb{R}^{N}}\Phi_1(|\nabla u|)dx
          + t^{l_{2}-\alpha}\int_{\mathbb{R}^{N}}\Phi_2(|\nabla v|)dx
          + t^{l_{1}-\alpha}\int_{\mathbb{R}^{N}}V_{1}(x)\Phi_1(| u|)dx
          \nonumber\\
&&\qquad\qquad
          + t^{l_{2}-\alpha}\int_{\mathbb{R}^{N}}V_{2}(x)\Phi_2(| v|)dx
          - t^{-\alpha}\int_{\mathbb{R}^{N}}\widetilde{F}(x,tu,tv)dx
          \nonumber\\
&&\qquad\quad
\geq       t^{l_{1}-\alpha}\int_{\mathbb{R}^{N}}\Phi_1(|\nabla u|)dx
          + t^{l_{2}-\alpha}\int_{\mathbb{R}^{N}}\Phi_2(|\nabla v|)dx
          + t^{l_{1}-\alpha}\int_{\mathbb{R}^{N}}V_{1}(x)\Phi_1(| u|)dx
          \nonumber\\
&&\qquad\qquad
          + t^{l_{2}-\alpha}\int_{\mathbb{R}^{N}}V_{2}(x)\Phi_2(| v|)dx
          - t^{-\alpha}C_4 \int_{\mathbb{R}^{N}}(t^{l_{1}r}|u|^{l_{1}r}+t^{l_{2}r}|v|^{l_{2}r})dx
          \nonumber\\
&&\qquad\quad
=           t^{l_{1}-\alpha}\int_{\mathbb{R}^{N}}\Phi_1(|\nabla u|)dx
          + t^{l_{2}-\alpha}\int_{\mathbb{R}^{N}}\Phi_2(|\nabla v|)dx
          + t^{l_{1}-\alpha}\int_{\mathbb{R}^{N}}V_{1}(x)\Phi_1(| u|)dx
          \nonumber\\
&&\qquad\qquad
          + t^{l_{2}-\alpha}\int_{\mathbb{R}^{N}}V_{2}(x)\Phi_2(| v|)dx
          - t^{l_{1}r-\alpha}C_4 \int_{\mathbb{R}^{N}}|u|^{l_{1}r}dx
          - t^{l_{2}r-\alpha}C_4 \int_{\mathbb{R}^{N}}|v|^{l_{2}r}dx
          \nonumber\\
&&\qquad\quad
\geq        t^{l_{1}-\alpha}
           \left(  \min\left\{\|\nabla u\|_{\Phi_1}^{l_{1}},\|\nabla u\|_{\Phi_1}^{m_{1}}\right\}
                 + \min\left\{\|u\|_{\Phi_1,V_{1}}^{l_{1}},\|u\|_{\Phi_1,V_{1}}^{m_{1}}\right\}
           \right)
          - t^{l_{1}r-\alpha}
            C_{5}
            \max\left\{\|u\|_{1,\Phi_1}^{l_{1}r},\|u\|_{1,\Phi_1}^{m_{1}r}\right\}
          \nonumber\\
&&\qquad\qquad
          + t^{l_{2}-\alpha}
           \left(  \min\left\{\|\nabla v\|_{\Phi_2}^{l_{2}},\|\nabla v\|_{\Phi_2}^{m_{2}}\right\}
                 + \min\left\{\|v\|_{\Phi_2,V_{2}}^{l_{2}},\|v\|_{\Phi_2,V_{2}}^{m_{2}}\right\}
           \right)
          - t^{l_{2}r-\alpha}
          C_{5}\max\left\{\|v\|_{1,\Phi_2}^{l_{2}r},\|v\|_{1,\Phi_2}^{m_{2}r}\right\}
          \nonumber\\
&&\qquad\quad
>0
\end{eqnarray}
because of $l_{1}>l_{1}r$ and $l_{2}>l_{2}r$ by the condition
$r\in \left(\max\left\{\frac{1}{l_{1}},\frac{1}{l_{2}}\right\},
\min\left\{\frac{l_{1}}{m_{1}},\frac{l_{2}}{m_{2}}\right\}\right)$.
Accordingly, for fixed $(u,v)\in W \backslash \{0\}$, $\gamma_{(u,v)}(t)$ has a unique zero $t(u,v)$ such that
$\gamma_{(u,v)}(t)<0$, for $0<t<t(u,v)$ and  $\gamma_{(u,v)}(t)\geq 0$, for $t(u,v)\leq t$.
Thus, the above result holds for $I_{0}(tu,tv)$, i.e., $I_{0}(u,v)$  satisfies condition $(A_{5})$. \qed
\par
\noindent
\begin{lemma}\label{lemma3.8}
Assume that $(\phi_1)$-$(\phi_4)$ and $(F_2')$ hold.  For any $b_{1}>0$ and $b_{2}>0$, there exist  $\sigma_{1}(b_{1})>0$ and $\sigma_{2}(b_{2})>0$ such that if $|\varepsilon|\leq \min\{\sigma_{1}(b_{1}),\sigma_{2}(b_{2})\}$, $I_{\varepsilon}'(u,v)=0$ and $|I_{\varepsilon}(u,v)|\leq \min\{\sigma_{1}(b_{1}),\sigma_{2}(b_{2})\}$, then $\|u\|_{W_{1}}\leq b_{1}$ and $\|v\|_{W_{2}}\leq b_{2}$.
\end{lemma}
\par
\noindent
{\bf Proof.}
Suppose on the contrary that there exist two sequences $\{(u_{n},v_{n})\}\subset W$ and $\{\varepsilon_{n}\}$ such that $\varepsilon_{n}\rightarrow 0$,
$I_{\varepsilon_{n}}(u_{n},v_{n})\rightarrow 0$ as $n\rightarrow\infty$, $I_{\varepsilon_{n}}'(u_{n},v_{n})=0$ and $\|(u_{n},v_{n})\|_{W}\geq b_{0}>0$, where $b_{0}$ is independent of $n$. Obviously, we can see $\{(u_{n},v_{n})\}$ as the (P-S) sequence of $I_{0}$. Then from Lemma \ref{lemma3.6}, we obtain that a subsequence of $\{(u_{n},v_{n})\}$ which converges to $(u_{0},v_{0})$ in $W$, satisfies
\begin{eqnarray*}\label{3.1.30}
&      \langle I_{0}'(u_{0},v_{0}),(u_{0},v_{0})\rangle
=     \int_{\mathbb{R}^{N}}\phi_1(|\nabla u_{0}|)\nabla u_{0} \nabla u_{0}dx
          + \int_{\mathbb{R}^{N}}\phi_2(|\nabla v_{0}|)\nabla v_{0} \nabla v_{0}dx
          + \int_{\mathbb{R}^{N}}V_{1}(x)\phi_1(|u_{0}|)u_{0}u_{0}dx
          \nonumber\\
&\qquad\qquad\qquad\qquad\qquad\quad
          + \int_{\mathbb{R}^{N}}V_{2}(x)\phi_2(|v_{0}|)v_{0}v_{0}dx
          - \int_{\mathbb{R}^{N}}\widetilde{F}_{u_{0}}(x,u_{0},v_{0})u_{0}dx
          - \int_{\mathbb{R}^{N}}\widetilde{F}_{v_{0}}(x,u_{0},v_{0})v_{0}dx
=0
\end{eqnarray*}
and
{\small\begin{eqnarray*}\label{3.1.31}
     I_{0}(u_{0},v_{0})
=     \int_{\mathbb{R}^{N}}\Phi_1(|\nabla u_{0}|)dx
          + \int_{\mathbb{R}^{N}}\Phi_2(|\nabla v_{0}|)dx
          + \int_{\mathbb{R}^{N}}V_{1}(x)\Phi_1(|u_{0}|)dx
          + \int_{\mathbb{R}^{N}}V_{2}(x)\Phi_2(|v_{0}|)dx
          - \int_{\mathbb{R}^{N}}\widetilde{F}(x,u_{0},v_{0})dx
=0.
\end{eqnarray*}}
From the above two equations, $(\phi_3)$ and $(F_2')$, it follows that
\begin{eqnarray}\label{3.1.32}
&&0=     I_{0}(u_{0},v_{0})-\frac{1}{\alpha}\langle I_{0}'(u_{0},v_{0}),(u_{0},v_{0})\rangle
        \nonumber\\
&&\quad
=     \int_{\mathbb{R}^{N}}\Phi_1(|\nabla u_{0}|)dx
    - \frac{1}{\alpha}\int_{\mathbb{R}^{N}}\phi_1(|\nabla u_{0}|)\nabla u_{0} \nabla u_{0}dx
    + \int_{\mathbb{R}^{N}}\Phi_2(|\nabla v_{0}|)dx
    - \frac{1}{\alpha}\int_{\mathbb{R}^{N}}\phi_2(|\nabla v_{0}|)\nabla v_{0} \nabla v_{0}dx
          \nonumber\\
&&\qquad
    + \int_{\mathbb{R}^{N}}V_{1}(x)\Phi_1(|u_{0}|)dx
    - \frac{1}{\alpha} \int_{\mathbb{R}^{N}}V_{1}(x)\phi_1(|u_{0}|)u_{0}u_{0}dx
    + \int_{\mathbb{R}^{N}}V_{2}(x)\Phi_2(|v_{0}|)dx
    - \frac{1}{\alpha}\int_{\mathbb{R}^{N}}V_{2}(x)\phi_2(|v_{0}|)v_{0}v_{0}dx
          \nonumber\\
&&\qquad
    + \frac{1}{\alpha}\int_{\mathbb{R}^{N}}\widetilde{F}_{u_{0}}(x,u_{0},v_{0})u_{0}dx
    + \frac{1}{\alpha}\int_{\mathbb{R}^{N}}\widetilde{F}_{v_{0}}(x,u_{0},v_{0})v_{0}dx
    - \int_{\mathbb{R}^{N}}\widetilde{F}(x,u_{0},v_{0})u_{0}dx
        \nonumber\\
&&\quad
\leq  \int_{\mathbb{R}^{N}}\Phi_1(|\nabla u_{0}|)dx
    - \frac{l_{1}}{\alpha}\int_{\mathbb{R}^{N}}\Phi_1(|\nabla u_{0}|)dx
    + \int_{\mathbb{R}^{N}}\Phi_2(|\nabla v_{0}|)dx
    - \frac{l_{2}}{\alpha}\int_{\mathbb{R}^{N}}\Phi_2(|\nabla v_{0}|)dx
          \nonumber\\
&&\qquad
    + \int_{\mathbb{R}^{N}}V_{1}(x)\Phi_1(|u_{0}|)dx
    - \frac{l_{1}}{\alpha} \int_{\mathbb{R}^{N}}V_{1}(x)\Phi_1(|u_{0}|)dx
    + \int_{\mathbb{R}^{N}}V_{2}(x)\Phi_2(|v_{0}|)dx
    - \frac{l_{2}}{\alpha}\int_{\mathbb{R}^{N}}V_{2}(x)\Phi_2(|v_{0}|)dx
          \nonumber\\
&&\qquad
    + \frac{C_{1,4}\left(l_{1}-\alpha\right)}{l_{1}\alpha}\int_{\mathbb{R}^{N}}V_{1}(x)|u_{0}|^{l_{1}}dx
    + \frac{C_{2,4}\left(l_{2}-\alpha\right)}{l_{2}\alpha}\int_{\mathbb{R}^{N}}V_{2}(x)|v_{0}|^{l_{2}}dx
        \nonumber\\
&&\quad
\leq  \int_{\mathbb{R}^{N}}\Phi_1(|\nabla u_{0}|)dx
    - \frac{l_{1}}{\alpha}\int_{\mathbb{R}^{N}}\Phi_1(|\nabla u_{0}|)dx
    + \int_{\mathbb{R}^{N}}\Phi_2(|\nabla v_{0}|)dx
    - \frac{l_{2}}{\alpha}\int_{\mathbb{R}^{N}}\Phi_2(|\nabla v_{0}|)dx
          \nonumber\\
&&\qquad
    - \frac{C_{1,4}\left(l_{1}-\alpha\right)}{\alpha}\int_{\mathbb{R}^{N}}V_{1}(x)|u_{0}|^{l_{1}}dx
    - \frac{C_{2,4}\left(l_{2}-\alpha\right)}{\alpha} \int_{\mathbb{R}^{N}}V_{2}(x)|v_{0}|^{l_{2}}dx
          \nonumber\\
&&\qquad
    + \frac{C_{1,4}\left(l_{1}-\alpha\right)}{l_{1}\alpha}\int_{\mathbb{R}^{N}}V_{1}(x)|u_{0}|^{l_{1}}dx
    + \frac{C_{2,4}\left(l_{2}-\alpha\right)}{l_{2}\alpha} \int_{\mathbb{R}^{N}}V_{2}(x)|v_{0}|^{l_{2}}dx.
\end{eqnarray}
Note that $l_{1}>\alpha$, $l_{2}>1$, $l_{2}>\alpha$, $l_{2}>1$ and $(u_{0},v_{0})$ in $W$.  The inequality \eqref{3.1.32} implies that $(u_{0},v_{0})\equiv(0,0)$.

On the other hand, from $\|(u_{n},v_{n})\|_{W}\geq b_{0}>0$ and $(u_{n},v_{n})\rightarrow (u_{0},v_{0})$ in $W$, we have
$\|(u_{0},v_{0})\|_{W}\geq b_{0}>0$, which contradicts the fact that $(u_{0},v_{0})\equiv (0,0)$. The proof is complete. \qed
\par
\noindent
\begin{lemma}\label{lemma3.9}
Assume that $(\phi_1)$-$(\phi_4)$ and $(F_1')$ hold.  There exist positive constants $C^{*}$, $D^{*}$, $\beta$ and $\eta$ such that  if $I_{\varepsilon}'(u,v)=0$ with $|\varepsilon|\leq 1$, then
$\|u\|_{\infty}\leq C^{*} \|u\|_{l_{1}^{\ast}}^{\beta}$ and $\|v\|_{\infty}\leq D^{*} \|v\|_{l_{2}^{\ast}}^{\eta}$.
\end{lemma}
\par
\noindent
{\bf Proof.}
Let $(u,v)\in W$ be a critical point of $I_{\varepsilon}$. Using  system \eqref{eq1}, for any $\psi=(\psi_{1},\psi_{2})\in  W$,
we obtain
\begin{eqnarray}\label{3.1.33}
&&           \int_{\mathbb{R}^{N}}\phi_1(|\nabla u|)\nabla u \nabla \psi_{1}dx
          + \int_{\mathbb{R}^{N}}\phi_2(|\nabla v|)\nabla v \nabla \psi_{2}dx
          + \int_{\mathbb{R}^{N}}V_{1}(x)\phi_1(|u|)u\psi_{1}dx
          + \int_{\mathbb{R}^{N}}V_{2}(x)\phi_2(|v|)v\psi_{2}dx\nonumber\\
&&
 =        \int_{\mathbb{R}^{N}}\widetilde{F}_{u}(x,u,v)\psi_{1}dx
          + \int_{\mathbb{R}^{N}}\widetilde{F}_{v}(x,u,v)\psi_{2}dx
          + \varepsilon\int_{\mathbb{R}^{N}}k(x)\widetilde{G}_{u}(u,v)\psi_{1}dx
          + \varepsilon\int_{\mathbb{R}^{N}}k(x)\widetilde{G}_{v}(u,v)\psi_{2}dx.
\end{eqnarray}
Choosing $\psi_{1}=|u^{T}|^{\nu_{1}}u^{T}$,  where  $\nu_{1}>0$, $u^{T}$ is defined as follows.
Take $T >\delta $ and define
\begin{align}\label{3.1.33-1}
      u^{T}
=
    \begin{cases}
         -T, \;\;\; \;\;\; \;\; \;             \text { if }\;\;u \leq -T, \\
          u, \;\;\; \;\;\; \;\; \;             \text { if }\;\;-T<u<T,\\
          T, \;\;\; \;\;\; \;\; \;             \text { if }\;\;u \geq T.
    \end{cases}
\end{align}
Choosing $\psi_{2}=|v^{T}|^{\nu_{2}}v^{T}$, where  $\nu_{2}>0$, $v^{T}$ is defined as follows.
Take $T >\delta $ and define
\begin{align}\label{3.1.33-2}
    v^{T}
=
    \begin{cases}
         -T, \;\;\; \;\;\; \;\; \;             \text { if }\;\;v \leq -T, \\
          v, \;\;\; \;\;\; \;\; \;             \text { if }\;\;-T<v<T,\\
          T, \;\;\; \;\;\; \;\; \;             \text { if }\;\;v \geq T.
    \end{cases}
\end{align}
By $(\phi_3)$ and \eqref{3.1.14},
\begin{align*}
           \int_{\mathbb{R}^{N}}\phi_1(|\nabla u|)\nabla u \nabla u(\nu_{1}+1)|u^{T}|^{\nu_{1}}dx
&\geq    (\nu_{1}+1)l_{1}          \int_{\mathbb{R}^{N}}\Phi_1(|\nabla u|)|u^{T}|^{\nu_{1}}dx
\geq    (\nu_{1}+1)l_{1}C_{1,4}   \int_{\mathbb{R}^{N}}|\nabla u|^{l_{1}}|u^{T}|^{\nu_{1}}dx
          \nonumber\\
&
\geq    (\nu_{1}+1)l_{1}C_{1,4}   \int_{\mathbb{R}^{N}}|\nabla u^{T}|^{l_{1}}|u^{T}|^{\nu_{1}}dx
=        \frac{(\nu_{1}+1)l_{1}^{l_{1}+1}C_{1,4}}{\left(\nu_{1}+l_{1}\right)^{l_{1}}}
         \int_{\mathbb{R}^{N}}\left|\nabla \left|u^{T}\right|^{\frac{\nu_{1}}{l_{1}}+1}\right|^{l_{1}}dx.
\end{align*}
Hence, by H\"{o}lder's inequality, $(F_1')$ and \eqref{3.1.14}, we have
\begin{eqnarray*}\label{3.1.37}
&&            \frac{(\nu_{1}+1)l_{1}^{l_{1}+1}C_{1,4}}{\left(\nu_{1}+l_{1}\right)^{l_{1}}}
         \int_{\mathbb{R}^{N}}\left|\nabla \left|u^{T}\right|^{\frac{\nu_{1}}{l_{1}}+1}\right|^{l_{1}}dx
          \nonumber\\
&&
\leq     C_{1}\int_{\mathbb{R}^{N}}|u|^{l_{1}r+\nu_{1}}dx
       + C_{7}\int_{\mathbb{R}^{N}}|u|^{1+\nu_{1}}dx
     + C_{1}
        \left(\int_{\mathbb{R}^{N}}|v|^{l_{2}}dx\right)^{\frac{l_{1}r-1}{l_{1}}}
        \left(\int_{\mathbb{R}^{N}}|u|^{\frac{l_{1}(1+\nu_{1})}{l_{1}(1-r)+1}}dx\right)^{\frac{l_{1}(1-r)+1}{l_{1}}}
          \nonumber\\
&&
\leq     C_{1}\int_{\mathbb{R}^{N}}|u|^{l_{1}r+\nu_{1}}dx
       + C_{7}\int_{\mathbb{R}^{N}}|u|^{1+\nu_{1}}dx
       + \frac{C_{1}}{C_{2,4}^{\frac{l_{1}r-1}{l_{1}}}}
        \left(\int_{\mathbb{R}^{N}}\Phi_{2}(|v|)dx\right)^{\frac{l_{1}r-1}{l_{1}}}
        \left(\int_{\mathbb{R}^{N}}|u|^{\frac{l_{1}(1+\nu_{1})}{l_{1}(1-r)+1}}dx\right)^{\frac{l_{1}(1-r)+1}{l_{1}}}
          \nonumber\\
&&
\leq     C_{1}\int_{\mathbb{R}^{N}}|u|^{l_{1}r+\nu_{1}}dx
       + C_{7}\int_{\mathbb{R}^{N}}|u|^{1+\nu_{1}}dx
     + \frac{C_{1}}{C_{2,4}^{\frac{l_{1}r-1}{l_{1}}}}
       \max\left\{\|v\|_{\Phi_{2}}^{\frac{l_{2}(l_{1}r-1)}{l_{1}}},\|v\|_{\Phi_{2}}^{\frac{m_{2}(l_{1}r-1)}{l_{1}}}\right\}
        \left(\int_{\mathbb{R}^{N}}|u|^{\frac{l_{1}(1+\nu_{1})}{l_{1}(1-r)+1}}dx\right)^{\frac{l_{1}(1-r)+1}{l_{1}}},
\end{eqnarray*}
where $C_{7}=\max_{|(u,v)|\leq \delta}|\widetilde{G}_{u}(u,v)|$.
Moreover, by Gagliard-Nirenberg-Sobolev inequality, we have
{\small\begin{eqnarray*}\label{3.1.38}
&&            \frac{(\nu_{1}+1)l_{1}^{l_{1}+1}C_{1,4}}{\left(\nu_{1}+l_{1}\right)^{l_{1}}}
         \left(\int_{\mathbb{R}^{N}}\left(|u^{T}|^{\frac{\nu_{1}}{l_{1}}+1}\right)^{l_{1}^{*}}dx\right)^{\frac{l_{1}}{l_{1}^{*}}}
\leq     \frac{(\nu_{1}+1)l_{1}^{l_{1}+1}C_{1,4}D_{1}}{\left(\nu_{1}+l_{1}\right)^{l_{1}}}
         \int_{\mathbb{R}^{N}}\left(|\nabla |u^{T}|^{\frac{\nu_{1}}{l_{1}}+1}|\right)^{l_{1}}dx
          \nonumber\\
&&
\leq     C_{1}D_{1}\int_{\mathbb{R}^{N}}|u|^{l_{1}r+\nu_{1}}dx
       + C_{7}D_{1}\int_{\mathbb{R}^{N}}|u|^{1+\nu_{1}}dx
       + \frac{C_{1}D_{1}}{C_{2,4}^{\frac{l_{1}r-1}{l_{1}}}}
         \max\left\{\|v\|_{\Phi_{2}}^{\frac{l_{2}(l_{1}r-1)}{l_{1}}},\|v\|_{\Phi_{2}}^{\frac{m_{2}(l_{1}r-1)}{l_{1}}}\right\}
        \left(\int_{\mathbb{R}^{N}}|u|^{\frac{l_{1}(1+\nu_{1})}{l_{1}(1-r)+1}}dx\right)^{\frac{l_{1}(1-r)+1}{l_{1}}}.
\end{eqnarray*}}
Next, we assume that
$\int_{\mathbb{R}^{N}}\left(|u|^{1+\nu_{1}}\right)^{q}dx=\max\left\{\int_{\mathbb{R}^{N}}|u|^{l_{1}r+\nu_{1}}dx,\int_{\mathbb{R}^{N}}|u|^{1+\nu_{1}}dx,
\int_{\mathbb{R}^{N}}|u|^{\frac{l_{1}(1+\nu_{1})}{l_{1}(1-r)+1}}dx\right\}$
, where $q=\frac{l_{1}}{l_{1}(1-r)+1}$ and $\frac{1}{q}=\frac{l_{1}(1-r)+1}{l_{1}}$,
the case $\int_{\mathbb{R}^{N}}|u|^{l_{1}r+\nu_{1}}dx=\max\left\{\int_{\mathbb{R}^{N}}|u|^{l_{1}r+\nu_{1}}dx,\int_{\mathbb{R}^{N}}|u|^{1+\nu_{1}}dx,
\int_{\mathbb{R}^{N}}|u|^{\frac{l_{1}(1+\nu_{1})}{l_{1}(1-r)+1}}dx\right\}$
and
$\int_{\mathbb{R}^{N}}|u|^{1+\nu_{1}}dx=\max\left\{\int_{\mathbb{R}^{N}}|u|^{l_{1}r+\nu_{1}}dx,\int_{\mathbb{R}^{N}}|u|^{1+\nu_{1}}dx,
\int_{\mathbb{R}^{N}}|u|^{\frac{l_{1}(1+\nu_{1})}{l_{1}(1-r)+1}}dx\right\}$ can be similarly treated.
So for fixed $u\in W_{1}$, we have the following estimate
\begin{eqnarray*}\label{3.1.39}
          \frac{(\nu_{1}+1)l_{1}^{l_{1}+1}C_{1,4}}{\left(\nu_{1}+l_{1}\right)^{l_{1}}}
         \left(\int_{\mathbb{R}^{N}}\left(|u^{T}|^{\frac{\nu_{1}}{l_{1}}+1}\right)^{l_{1}^{*}}dx\right)^{\frac{l_{1}}{l_{1}^{*}}}
\leq     D_{2}\max\left\{\|v\|_{\Phi_{2}}^{\frac{l_{2}(l_{1}r-1)}{l_{1}}},\|v\|_{\Phi_{2}}^{\frac{m_{2}(l_{1}r-1)}{l_{1}}}\right\}
         \left(\int_{\mathbb{R}^{N}}\left(|u|^{1+\nu_{1}}\right)^{q}dx\right)^{\frac{1}{q}},
\end{eqnarray*}
where $D_{2}:=3D_{1}\max\left\{C_{1},C_{7},\frac{C_{1}}{C_{2,4}^{\frac{l_{1}r-1}{l_{1}}}} \right\}$.
Taking $T\rightarrow\infty$, it follows that
\begin{eqnarray*}\label{3.1.40}
\frac{l_{1}^{l_{1}+1}C_{1,4}}{\left(\nu_{1}+l_{1}\right)^{l_{1}-1}}
\left(\int_{\mathbb{R}^{N}}\left(|u|^{\frac{\nu_{1}}{l_{1}}+1}\right)^{l_{1}^{*}}dx\right)^{\frac{l_{1}}{l_{1}^{*}}}
\leq     D_{2}
\max\left\{\|v\|_{\Phi_{2}}^{\frac{l_{2}(l_{1}r-1)}{l_{1}}},\|v\|_{\Phi_{2}}^{\frac{m_{2}(l_{1}r-1)}{l_{1}}}\right\}
         \left(\int_{\mathbb{R}^{N}}\left(|u|^{1+\nu_{1}}\right)^{q}dx\right)^{\frac{1}{q}},
\end{eqnarray*}
i.e.,
\begin{eqnarray*}\label{3.1.41}
          \frac{l_{1}^{l_{1}+1}C_{1,4}}{\left(\nu_{1}+l_{1}\right)^{l_{1}-1}}
          \|u\|_{\frac{(\nu_{1}+l_{1})N}{N-l_{1}}}^{\nu_{1}+l_{1}}
\leq      D_{2}
\max\left\{\|v\|_{\Phi_{2}}^{\frac{l_{2}(l_{1}r-1)}{l_{1}}},\|v\|_{\Phi_{2}}^{\frac{m_{2}(l_{1}r-1)}{l_{1}}}\right\}
          \|u\|_{(1+\nu_{1})q}^{1+\nu_{1}}.
\end{eqnarray*}
So,
\begin{eqnarray*}\label{3.1.42}
&&          \|u\|_{\frac{(\nu_{1}+l_{1})N}{N-l_{1}}}
\leq    \left(\nu_{1}+l_{1}\right)^{\frac{l_{1}}{\nu_{1}+l_{1}}}
        \left(\frac{D_{2}}{l_{1}^{l_{1}+1}C_{1,4}}\right)^{\frac{1}{\nu_{1}+l_{1}}}
        \max\left\{\|v\|_{\Phi_{2}}^{\frac{l_{2}(l_{1}r-1)}{l_{1}(\nu_{1}+l_{1})}},\|v\|_{\Phi_{2}}^{\frac{m_{2}(l_{1}r-1)}{l_{1}(\nu_{1}+l_{1})}}\right\}
        \|u\|_{(1+\nu_{1})q}^{\frac{1+\nu_{1}}{\nu_{1}+l_{1}}}\nonumber\\
&&\qquad\qquad
\leq    \left(\nu_{1}+l_{1}\right)^{\frac{l_{1}}{\nu_{1}+l_{1}}}D_{3}^{\frac{l_{1}}{\nu_{1}+l_{1}}}
        \|u\|_{(1+\nu_{1})q}^{\frac{1+\nu_{1}}{\nu_{1}+l_{1}}}
        =\left(D_{3}(\nu_{1}+l_{1})\right)^{\frac{l_{1}}{\nu_{1}+l_{1}}}\|u\|_{(1+\nu_{1})q}^{\frac{1+\nu_{1}}{\nu_{1}+l_{1}}},
\end{eqnarray*}
where
$D_{3}=
\max\left\{1,\left(\frac{D_{2}}{l_{1}^{l_{1}+1}C_{1,4}}\right)^{\frac{1}{\nu_{1}+l_{1}}}
\max\left\{\|v\|_{\Phi_{2}}^{\frac{l_{2}(l_{1}r-1)}{l_{1}(\nu_{1}+l_{1})}},
\|v\|_{\Phi_{2}}^{\frac{m_{2}(l_{1}r-1)}{l_{1}(\nu_{1}+l_{1})}}\right\}\right\}$.
So $D_{3}\geq 1$.
Let us define $\nu_{1,k}=\frac{(\nu_{1,k-1}+l_{1})N}{q(N-l_{1})}-1$, where $k=1,2,\cdots$ and $\nu_{1,0}=\frac{l_{1}^{\ast}-q}{q}$, that is
$\nu_{1,k}=\frac{\left(\frac{l_{1}^{\ast}}{ql_{1}}\right)^{k+1}-1}{\left(\frac{l_{1}^{\ast}}{ql_{1}}\right)-1}\nu_{1,0}$,
for $k=1,2,\cdots$, (see  Lemma \ref{lemma6.1}  in Appendix).
Let $a=\frac{l_{1}^{\ast}}{ql_{1}}$. Then $a>1$ by $l_{i}^{\ast}m_i(l_i+1)-l_{i}^{2}(l_{i}^{\ast}+m_i)\geq 0$ (see Lemma \ref{lemma6.2}  in Appendix).
It is easy to see that $\nu_{1,k}\rightarrow +\infty$ as $k\rightarrow +\infty$.

By Moser's iteration method we have
{\footnotesize\begin{align*}
\|u\|_{\frac{(\nu_{1,k}+l_{1})N}{N-l_{1}}}
&\leq    \left(D_{3}(\nu_{1,k}+l_{1})\right)^{\frac{l_{1}}{\nu_{1,k}+l_{1}}}\cdots
        \left(D_{3}(\nu_{1,1}+l_{1})\right)^{\frac{l_{1}}{\nu_{1,1}+l_{1}}}
        \left(D_{3}(\nu_{1,0}+l_{1})\right)^{\frac{l_{1}}{\nu_{1,0}+l_{1}}}
        \|u\|_{l_{1}^{\ast}}^{\frac{1+\nu_{1,0}}{\nu_{1,0}+l_{1}}\frac{1+\nu_{1,1}}{\nu_{1,1}+l_{1}}\frac{1+\nu_{1,2}}{\nu_{1,2}+l_{1}}\cdots\frac{1+\nu_{1,k}}{\nu_{1,k}+l_{1}}}
        \nonumber\\
&
\leq    \exp\left(\sum_{i=0}^{k}\frac{l_{1}\ln\left(D_{2}(\nu_{1,i}+l_{1}\right)}{\nu_{1,i}+l_{1}}\right)
        \|u\|_{l_{1}^{\ast}}^{\beta_{k}},
\end{align*}}
where $\beta_{k}=\Pi_{i=0}^{k}\frac{1+\nu_{1,i}}{\nu_{1,i}+l_{1}}$.
Taking $k\rightarrow\infty$, we obtain that
\begin{eqnarray*}
         \|u\|_{\infty}
\leq
       C^{*} \|u\|_{l_{1}^{\ast}}^{\beta},
\end{eqnarray*}
where $0<\beta=\Pi_{i=0}^{\infty}\frac{1+\nu_{1,i}}{\nu_{1,i}+l_{1}}<1$, $C^{*}=\exp\left(\sum_{i=0}^{\infty}\frac{l_{1}\ln\left(D_{2}(\nu_{1,i}+l_{1}\right)}{\nu_{1,i}+l_{1}}\right)$ which is a positive constant, (for details, see  Lemma \ref{lemma6.3} and Lemma \ref{lemma6.4} in appendix).
Similarly, we can get
\begin{eqnarray*}
         \|v\|_{\infty}
\leq
       D^{*} \|v\|_{l_{2}^{\ast}}^{\eta}.
\end{eqnarray*}
The proof is complete.\qed
\par
\noindent
\begin{remark}\label{remark3.10}
By Lemma \ref{lemma3.8}, if $(u,v)$ is critical point of $I_{\varepsilon}(u,v)$ with
$|\varepsilon|\leq \min\left\{\sigma_{1}\left(\left(\frac{\sqrt{2}\delta}{4C^{*}C_{1,3}^{\beta}}\right)^{\frac{1}{\beta}}\right),
\sigma_{2}\left(\left(\frac{\sqrt{2}\delta}{4D^{*}C_{2,3}^{\eta}}\right)^{\frac{1}{\eta}}\right)\right\}$
and
$|I_{\varepsilon}(u,v)|\leq
\min\left\{\sigma_{1}\left(\left(\frac{\sqrt{2}\delta}{4C^{*}C_{1,3}^{\beta}}\right)^{\frac{1}{\beta}}\right),
\sigma_{2}\left(\left(\frac{\sqrt{2}\delta}{4D^{*}C_{2,3}^{\eta}}\right)^{\frac{1}{\eta}}\right)\right\}$,
then $\|u\|_{W_{1}}\leq \left(\frac{\sqrt{2}\delta}{4C^{*}C_{1,3}^{\beta}}\right)^{\frac{1}{\beta}}$
and
$\|v\|_{W_{2}}\leq  \left(\frac{\sqrt{2}\delta}{4D^{*}C_{2,3}^{\eta}}\right)^{\frac{1}{\eta}}$
which together with Lemma \ref{lemma3.9} implies that
$\|u\|_{\infty}\leq \frac{\sqrt{2}\delta}{4}$
and
$\|v\|_{\infty}\leq \frac{\sqrt{2}\delta}{4}$.
This means $(u,v)$ is critical point of the original problem  \eqref{eq1}.
\end{remark}
\par
\vskip2mm
\noindent
{\bf  Proof of Theorem 1.1.}
Without loss of generality, we assume that $\varepsilon>0$. The case
$\varepsilon<0$ can be studied similarly by replacing $\widetilde{G}_u(u,v)$ and $\widetilde{G}_v(u,v)$
with $-\widetilde{G}_u(u,v)$ and $-\widetilde{G}_v(u,v)$.
\par
Next we are ready to verify that $I_{\varepsilon}(u,v)$ satisfies conditions $(A_{1})$-$(A_{5})$ in Lemma \ref{lemma2.10}.
To verify condition $(A_{1})$, by \eqref{3.1.23}, we have
\begin{eqnarray*}
 \inf_{\varepsilon\in [0,1],(u,v)\in W}I_{\varepsilon}(u,v)>-\infty.
\end{eqnarray*}
Condition $(A_{2})$ follows from
\begin{eqnarray*}
       |I_{\varepsilon}(u,v)-I_{0}(u,v)|
\leq   |\varepsilon|\int_{\mathbb{R}^{N}}|k(x)||\widetilde{G}(u,v)|dx
\leq   |\varepsilon|C:=\psi(\varepsilon),
\end{eqnarray*}
where $C$ is a constant independent of $(u,v)$ and $\varepsilon$.
Conditions $(A_{3})$ and $(A_{5})$ follow from Lemma \ref{lemma3.6} and Lemma \ref{lemma3.7}, respectively.
Thus $I_{\varepsilon}(u,v)$ satisfies all the conditions in Lemma \ref{lemma2.10}. Then by using the same proofs as Corollary 1.1 in \cite{Huang2022}, for any $\min\{\sigma_{1},\sigma_{2}\}>0$ and any given $k\in \mathbb N$  , we have $k$ distinct critical values of $I_{\varepsilon}$ satisfying
\begin{eqnarray*}
       -\min\{\sigma_{1},\sigma_{2}\}
<   a_{n(1)}(\varepsilon)
<   a_{n(2)}(\varepsilon)
<
\cdots
<   a_{n(k)}(\varepsilon)
<0.
\end{eqnarray*}
Finally, due to the arbitrariness of $\min\{\sigma_{1},\sigma_{2}\}$,
we can take $0<\sigma<\min\left\{\sigma_{1}\left(\left(\frac{\sqrt{2}\delta}{4C^{*}C_{1,3}^{\beta}}\right)^{\frac{1}{\beta}}\right),
\sigma_{2}\left(\left(\frac{\sqrt{2}\delta}{4D^{*}C_{2,3}^{\eta}}\right)^{\frac{1}{\beta}}\right)\right\}$.
Then by Remark \ref{remark3.10}, the original problem  \eqref{eq1} has at least $k$ solutions whose $L^{\infty}$-norms are less that $\frac{\sqrt{2}\delta}{2}$. \qed
\par
\noindent
\section{Example}

\begin{eqnarray}\label{eq6}
 \begin{cases}
-\mbox{div}[(4|\nabla u|^{2}+5|\nabla u|^{3})\nabla u]+V_{1}(x)(4|u|^{2}+5|u|^{3})u= F_u(x,u,v), \ \ x\in \mathbb{R}^6,\\
-\mbox{div}\left[\left(4|\nabla v|^{2}\ln(e+|\nabla v|)+\dfrac{|\nabla v|^{3}}{e+|\nabla v|}\right)\nabla v\right]
+V_{1}(x)\left(4|v|^{2}\ln(e+| v|)+\dfrac{|v|^{3}}{e+|v|}\right)v
= F_v(x,u,v), \ \ x\in \mathbb{R}^6,\\
 \end{cases}
\end{eqnarray}
where
\begin{eqnarray}\label{4.1.1}
F(x,t,s)= (\sin x_{1}+2)(|t|^{3}+|s|^{3}+t^{6}s^{6}).
\end{eqnarray}
Let $N=6$, $\phi_1(t)=4|t|^{2}+5|t|^{3}$ and $\phi_2(t)=4|t|^{2}\ln(e+| t|)+\dfrac{|t|^{3}}{e+|t|}$. Then $\phi_i (i=1,2)$ satisfy $(\phi_1)$-$(\phi_4)$, $l_{1}=l_{2}=4$, $m_{1}=5$, $m_{2}=4+\frac{1}{1+eln(e-1)}$ and $\Phi_1(t)=|t|^{4}+|t|^{5}$, $\Phi_2(t)=|t|^{4}\ln(e+|t|)$. So, $l^{*}_{1}=l^{*}_{2}=12$ and $4<m_{2}<5$.
\par
 Let $V_1(x)=\sum_{i=1}^{6}x_i^{2}+75$ and $V_2(x)=\sum_{i=1}^{6}x_i^{4}+150$
 for all $(x,t,s)\in \mathbb{R}^{N}\times \mathbb{R}\times \mathbb{R}$. Then it is obvious that $V_i, i=1,2$ satisfy (V0) and (V1).
 \par
By (\ref{4.1.1}), we have
\begin{eqnarray*}
&F_{t}(x,t,s)=(\sin x_{1}+2)\left(3|t|t+6t^{5}s^{6}\right), \label{4.1.2}\\
&F_{s}(x,t,s)=(\sin x_{1}+2)\left(3|s|s+6s^{5}t^{6}\right).\label{4.1.3}
\end{eqnarray*}
Hence
\begin{eqnarray*}
&|F_{t}(x,t,s)|\leq 9(t^{2}+ 2s^{2}), \;\;\text{for\;all}\;\;|(t,s)|<1, \label{4.1.4}\\
&|F_{s}(x,t,s)|\leq 9(s^{2}+ 2t^{2}), \;\;\text{for\;all}\;\;|(t,s)|<1,\label{4.1.5}
\end{eqnarray*}
where we choose $r=\frac{3}{4}$.
Besides,
\begin{align*}
&      F_{t}(x,t,s)t+F_{s}(x,t,s)s-\alpha F(x,t,s)\\
&
=      (\sin x_{1}+2)\left(3|t|^{3}+6t^{6}s^{6}\right)+(\sin x_{1}+2)\left(3|s|^{3}+6t^{6}s^{6}\right)
     - \alpha (\sin x_{1}+2)\left(|t|^{3}+|s|^{3}+t^{6}s^{6}\right)\\
&
=     (\sin x_{1}+2)(3-\alpha)(|t|^{3}+|s|^{3})+ (\sin x_{1}+2)(12-\alpha)t^{6}s^{6}
     \\
&
\leq   25t^{6}s^{6}
\leq   \frac{25}{2}t^{12}+\frac{25}{2}s^{12}
\leq      \frac{25}{2}|t|^{4}+\frac{25}{2}|s|^{4}
\\
&
\leq
\frac{\Phi_1(1)\left(l_{1}-\alpha\right)}{l_{1}}V_{1}(x)|t|^{4}
    +\frac{\Phi_2(1)\left(l_{2}-\alpha\right)}{l_{2}} V_{2}(x)|s|^{4}, \;\;\text{for\;all}\;\;|(t,s)|<1,
\end{align*}
where $\alpha=\frac{11}{3}$.
It is easy to see that conditions $(F_3)$ and $(F_4)$ hold.
Since
\begin{align*}
      \lim_{|(t,s)|\rightarrow 0}
      \left(\inf_{x\in \mathbb{R}^{N}}\frac{F(x,t,s)}{|t|^{\frac{11}{3}}+|s|^{\frac{11}{3}}}\right)
&=     \lim_{|(t,s)|\rightarrow 0}
      \left(\inf_{x\in \mathbb{R}^{N}}
      \frac{(\sin x_{1}+2)(|t|^{3}+|s|^{3}+t^{6}s^{6})}{|t|^{\frac{11}{3}}+|s|^{\frac{11}{3}}}\right)
=     \lim_{|(t,s)|\rightarrow 0}
      \frac{|t|^{3}+|s|^{3}+t^{6}s^{6}}{|t|^{\frac{11}{3}}+|s|^{\frac{11}{3}}}\\
&\leq  \lim_{|(t,s)|\rightarrow 0}\left(|t|^{-\frac{2}{3}}+|s|^{-\frac{2}{3}}\right)
      +\lim_{|(t,s)|\rightarrow 0}\frac{t^{6}s^{6}}{|t|^{\frac{11}{3}}+|s|^{\frac{11}{3}}}
=\infty.
\end{align*}
So the condition $(F_4)$ holds. Then by Corollary 1.2, system (\ref{eq6}) has infinitely many nontrivial weak solutions.
\section{Results for the scalar equation}
In this section, we study the multiplicity of solutions for the following perturbation elliptic equation in Orlicz-Sobolev spaces:
 \begin{equation}\label{eq1-1}
 \left\{
  \begin{array}{ll}
-\Delta_{\Phi}u+V(x)\phi(|u|)u=f(x,u)+\varepsilon k(x)g(u), &x\in \mathbb{R}^N,\\
  u\in W^{1,\Phi}(\mathbb{R}^N),
    \end{array}
 \right.
 \end{equation}
where  $\Delta_{\Phi}(u)=\mbox{div}(\phi(|\nabla u|)\nabla u)$, $\Phi(t):=\int_{0}^{|t|}s\phi(s)ds, \ t\in \mathbb{R}$, $\phi  :(0,+\infty)\rightarrow(0,+\infty)$ is a function
satisfying
\begin{itemize}
\item[$(\phi_1)'$] $\phi\in C^1(0,+\infty)$, $t\phi(t)\rightarrow0$ as
 $t\rightarrow0$, $t\phi(t)\rightarrow+\infty$ as
 $t\rightarrow+\infty$;
 \item[$(\phi_2)'$] $t\rightarrow t\phi(t)$ are strictly increasing;
 \item[$(\phi_3)'$] $1<l:=\inf_{t>0}\frac{t^2\phi(t)}{\Phi(t)}\leq\sup_{t>0}\frac{t^2\phi(t)}{\Phi(t)}=:m<\min\{N,l^{\ast},l^{2}\}$,
      where $l^{\ast}=\frac{lN}{N-l}$ and $l^{\ast}m(l+1)-l^{2}(l^{\ast}+m)\geq 0$;
 \item[$(\phi_4)'$]
    there exist positive constants $c_{1,1}$ and $c_{1,2}$  such that
    $$
    c_{1,1}|t|^{l}\le \Phi(t)\le c_{1,2}|t|^{l},\ \ \forall \;|t|<1.
    $$
\end{itemize}
\par
Moreover, to state our result, we introduce the following  assumptions concerning $f$, $V$, $k$ and $g$:
\begin{itemize}
\item[$(f_0)$]  $f: \mathbb{R}^N\times [-\delta, \delta]\rightarrow \mathbb{R}$ is a continuous;
\item[$(f_1)$] there exists a  constant $c_1>0$ such that
\begin{equation*}\label{3.1.1}
 |f(x,u)|\leq c_1|u|^{l_{1}r-1}
 \end{equation*}
for all $|u|<\delta$  and $x\in  \mathbb{R}^{N}$, where
$r\in \left[\frac{1}{l}, \frac{l}{m}\right)$;
 \item[$(f_2)$] there exists  a constant $\alpha\in [lr,l)$ such that
\begin{eqnarray*}
      f(x,u)u-\alpha F(x,u)
\leq a(x)|u|^{l}
\end{eqnarray*}
for all $|u|<\delta$  and $x\in  \mathbb{R}^{N}$, where $a(x)=\frac{\Phi(1)\left(l-\alpha\right)}{l}V(x)$, $F(x,u)=\int_{0}^{u}f(x,s)ds$;
\item[$(f_3)$] $\lim_{|u|\rightarrow 0}\left(\inf_{x \in \mathbb{R}^{N}}\frac{F(x,u)}{|u|^{\alpha}}\right)=\infty$;
\item[$(f_4)$]
 $f(x,-u)=-f(x,u)$
for all $|u|<\delta$  and $x\in  \mathbb{R}^{N}$;
\end{itemize}
\begin{itemize}
\item[$(V)'$] $V \in C(\mathbb{R}^{N},\mathbb{R})$,
$V(x)^{-1}\in L^{\frac{r}{1-r}}(\mathbb{R}^{N})$
\mbox{and} there exists a constant $c_{2}$ such that
$0<c_{2}<V(x)$;
\end{itemize}
\begin{itemize}
\item[$(K)'$] $k(x)\in L^{1}(\mathbb{R}^N)\cap L^{\infty}(\mathbb{R}^N)$.
\end{itemize}
\begin{itemize}
\item[$(G)'$] $g: [-\delta, \delta]\rightarrow \mathbb{R}$ is a continuous function.
\end{itemize}
\par
By similar proofs of Theorem 1.1, we can obtain the following result.
\par
\begin{theorem}\label{theorem5.1}
  Assume that $(\phi_1)$--$(\phi_4)$, $(f_0)$--$(f_4)$, $(V)'$, $(G)'$ and $(K)'$ hold.
 Then
 for any $ k \in \mathbb{N}$ and  any $b>0$, there exists an  $\sigma_{0}(k,b)>0$ such that when $|\varepsilon|\leq \sigma_{0}(k,b)$,  equation \eqref{eq1-1} possesses at least $k$ distinct solutions whose
 $L^{\infty}$-norms are less than $\frac{\sqrt{2}\delta}{4}$.
 \end{theorem}

\begin{corollary}
 {\it Assume that $(\phi_1)$--$(\phi_4)$, $(f_0)$--$(f_4)$ and $(V)'$ hold.
 Then equation \eqref{eq1-1} with $\varepsilon=0$ possesses infinitely many distinct solutions whose
 $L^{\infty}$-norms are less than $\frac{\sqrt{2}\delta}{4}$.}
 \end{corollary}

\section{Appendix}
\begin{lemma}\label{lemma6.1}
If $\nu_{1,k}= \frac{(\nu_{1,k-1}+l_{1})N}{q(N-l_{1})}-1$,
where $k=1,2,\cdots$ and $\nu_{1,0}=\frac{l_{1}^{\ast}-q}{q}=\frac{l_{1}^{\ast}}{q}-1$, then
$\nu_{1,k}=\frac{\left(\frac{l_{1}^{\ast}}{ql_{1}}\right)^{k+1}-1}{\left(\frac{l_{1}^{\ast}}{ql_{1}}\right)-1}\nu_{1,0}$,
for $k=1,2,\cdots$.
\end{lemma}
\par
\noindent
{\bf Proof.}
Since $\nu_{1,k}=\frac{1}{q}\left[\frac{(\nu_{1,k-1}+l_{1})N}{N-l_{1}}-q\right]$,
where $k=1,2,\cdots$ and $\nu_{1,0}=\frac{l_{1}^{\ast}-q}{q}=\frac{l_{1}^{\ast}}{q}-1$,
it is easy to see that
\begin{align*}
k=1 \Rightarrow
  \nu_{1,1}
&
= \frac{1}{q}\left[\frac{(\nu_{1,0}+l_{1})N}{N-l_{1}}-q\right]
= \frac{1}{q}\left[\frac{(\nu_{1,0}+l_{1})l_{1}^{\ast}}{l_{1}}-q\right]
= \frac{l_{1}^{\ast}\nu_{1,0}}{ql_{1}}+\frac{l_{1}^{\ast}}{q}-1
= \frac{l_{1}^{\ast}\nu_{1,0}}{ql_{1}}+\nu_{1,0}
  \nonumber\\
&
= \left(\frac{l_{1}^{\ast}}{ql_{1}}+1\right)\nu_{1,0}
= \frac{\left(\frac{l_{1}^{\ast}}{ql_{1}}\right)^{2}-1}{\frac{l_{1}^{\ast}}{ql_{1}}-1}\nu_{1,0},
\end{align*}
\begin{align*}
k=2 \Rightarrow
  \nu_{1,2}
&
= \frac{1}{q}\left[\frac{(\nu_{1,1}+l_{1})N}{N-l_{1}}-q\right]
= \frac{1}{q}\left[\frac{(\nu_{1,1}+l_{1})l_{1}^{\ast}}{l_{1}}-q\right]
= \frac{l_{1}^{\ast}\nu_{1,1}}{ql_{1}}+\frac{l_{1}^{\ast}}{q}-1
= \frac{l_{1}^{\ast}}{ql_{1}}\nu_{1,1}+\nu_{1,0}
  \nonumber\\
&
= \frac{l_{1}^{\ast}}{ql_{1}}\left(\frac{l_{1}^{\ast}}{ql_{1}}+1\right)\nu_{1,0}+\nu_{1,0}
= \left(\left(\frac{l_{1}^{\ast}}{ql_{1}}\right)^{2}+\frac{l_{1}^{\ast}}{ql_{1}}+1\right)\nu_{1,0}
= \frac{\left(\frac{l_{1}^{\ast}}{ql_{1}}\right)^{3}-1}{\frac{l_{1}^{\ast}}{ql_{1}}-1}\nu_{1,0},
\end{align*}
\begin{align*}
k=3 \Rightarrow
  \nu_{1,3}
&
= \frac{1}{q}\left[\frac{(\nu_{1,2}+l_{1})N}{N-l_{1}}-q\right]
= \frac{1}{q}\left[\frac{(\nu_{1,2}+l_{1})l_{1}^{\ast}}{l_{1}}-q\right]
= \frac{l_{1}^{\ast}\nu_{1,2}}{ql_{1}}+\frac{l_{1}^{\ast}}{q}-1
= \frac{l_{1}^{\ast}}{ql_{1}}\nu_{1,2}+\nu_{1,0}
  \nonumber\\
&
= \frac{l_{1}^{\ast}}{ql_{1}}
  \left(\left(\frac{l_{1}^{\ast}}{ql_{1}}\right)^{2}+\frac{l_{1}^{\ast}}{ql_{1}}+1\right)\nu_{1,0}+\nu_{1,0}
= \left( \left(\frac{l_{1}^{\ast}}{ql_{1}}\right)^{3}+\left(\frac{l_{1}^{\ast}}{ql_{1}}\right)^{2}
        +\frac{l_{1}^{\ast}}{ql_{1}}+1
  \right)\nu_{1,0}
= \frac{\left(\frac{l_{1}^{\ast}}{ql_{1}}\right)^{4}-1}{\frac{l_{1}^{\ast}}{ql_{1}}-1}\nu_{1,0}.
\end{align*}
Then by applying the mathematical induction, we have
\begin{align*}
  \nu_{1,k}
&
= \frac{1}{q}\left[\frac{(\nu_{1,k}+l_{1})N}{N-l_{1}}-q\right]
= \frac{1}{q}\left[\frac{(\nu_{1,k}+l_{1})l_{1}^{\ast}}{l_{1}}-q\right]
= \frac{l_{1}^{\ast}\nu_{1,k}}{ql_{1}}+\frac{l_{1}^{\ast}}{q}-1
= \frac{l_{1}^{\ast}}{ql_{1}}\nu_{1,k}+\nu_{1,0}
  \nonumber\\
&
= \frac{l_{1}^{\ast}}{ql_{1}}
  \left(\left(\frac{l_{1}^{\ast}}{ql_{1}}\right)^{k-1}+\cdots+\frac{l_{1}^{\ast}}{ql_{1}}+1\right)\nu_{1,0}+\nu_{1,0}
= \left( \left(\frac{l_{1}^{\ast}}{ql_{1}}\right)^{k}+\cdots+\left(\frac{l_{1}^{\ast}}{ql_{1}}\right)^{2}
        +\frac{l_{1}^{\ast}}{ql_{1}}+1
  \right)\nu_{1,0}
= \frac{\left(\frac{l_{1}^{\ast}}{ql_{1}}\right)^{k+1}-1}{\frac{l_{1}^{\ast}}{ql_{1}}-1}\nu_{1,0}
\end{align*}
The proof is complete. \qed
\begin{lemma}\label{lemma6.2}
If $q=\frac{l_{1}}{l_{1}(1-r)+1}$, $l_{1}^{\ast}m_1(l_1+1)-l_{1}^{2}(l_{1}^{\ast}+m_1)\geq 0$ and
$r\in \left(\max\left\{\frac{1}{l_{1}},\frac{1}{l_{2}}\right\}, \min\left\{\frac{l_{1}}{m_{1}},\frac{l_{2}}{m_{2}}\right\}\right)$,
then
$\frac{l_{1}^{\ast}}{ql_{1}}>1$.
\end{lemma}
\par
\noindent
{\bf Proof.}
Since $\frac{1}{q}=\frac{l_{1}(1-r)+1}{l_{1}}$ and
$r\in \left(\max\left\{\frac{1}{l_{1}},\frac{1}{l_{2}}\right\}, \min\left\{\frac{l_{1}}{m_{1}},\frac{l_{2}}{m_{2}}\right\}\right)$,
it is easy to see that
\begin{align}\label{6.1.1}
\frac{l_{1}^{\ast}}{ql_{1}}
\in
\left(
\frac{l_{1}^{\ast}\left(l_{1}-l_{1}\min\left\{\frac{l_{1}}{m_{1}},\frac{l_{2}}{m_{2}}\right\}+1\right)}{l_{1}^{2}},
\frac{l_{1}^{\ast}\left(l_{1}-l_{1}\max\left\{\frac{1}{l_{1}},\frac{1}{l_{2}}\right\}+1\right)}{l_{1}^{2}}
\right).
\end{align}
Since
$l_{1}^{\ast}m_1(l_1+1)-l_{1}^{2}(l_{1}^{\ast}+m_1)\geq 0$,
 we have
\begin{align*}
&
\frac{l_{1}^{\ast}m_1l_i+l_{1}^{\ast}m_1-l_{1}^{2}l_{1}^{\ast}+l_{1}^{2}m_1}{l_{1}^{\ast}l_1m_1}\geq 0
\Rightarrow
1+\frac{1}{l_{1}}-\frac{l_{1}}{m_{1}}-\frac{l_{1}}{l_{1}^{\ast}}\geq 0
\Rightarrow
1-\frac{l_{1}}{l_{1}^{\ast}}+\frac{1}{l_{1}}-\min\left\{\frac{l_{1}}{m_{1}},\frac{l_{2}}{m_{2}}\right\}\geq 0
\Rightarrow
  \nonumber\\
&
1+\frac{1}{l_{1}}-\min\left\{\frac{l_{1}}{m_{1}},\frac{l_{2}}{m_{2}}\right\}\geq \frac{l_{1}}{l_{1}^{\ast}}
\Rightarrow
\frac{l_{1}-l_{1}\min\left\{\frac{l_{1}}{m_{1}},\frac{l_{2}}{m_{2}}\right\}+1}{l_{1}}\geq \frac{l_{1}}{l_{1}^{\ast}}
\Rightarrow
l_{1}^{\ast}\left(l_{1}-l_{1}\min\left\{\frac{l_{1}}{m_{1}},\frac{l_{2}}{m_{2}}\right\}+1\right)\geq l_{1}^{2}
  \nonumber\\
&
\Rightarrow
\frac{l_{1}^{\ast}\left(l_{1}-l_{1}\min\left\{\frac{l_{1}}{m_{1}},\frac{l_{2}}{m_{2}}\right\}+1\right)}{l_{1}^{2}}\geq1.
\end{align*}
Hence,
let $a=\frac{l_{1}^{\ast}}{ql_{1}}$. Then $a>1$ by \eqref{6.1.1} and $l_{i}^{\ast}m_i(l_i+1)-l_{i}^{2}(l_{i}^{\ast}+m_i)\geq 0$. \qed

\begin{lemma}\label{lemma6.3}
The positive series $\sum_{i=0}^{\infty}\frac{l_{1}\ln\left(D_{2}(\nu_{1,i}+l_{1})\right)}{\nu_{1,i}+l_{1}}$
is convergent.
\end{lemma}
\par
\noindent
{\bf Proof.} Note that
\begin{align*}
   \sum_{i=0}^{\infty}\frac{l_{1}\ln\left(D_{2}(\nu_{1,i}+l_{1})\right)}{\nu_{1,i}+l_{1}}
=  \sum_{i=1}^{\infty}\frac{l_{1}\ln\left(D_{2}(\nu_{1,i}+l_{1})\right)}{\nu_{1,i}+l_{1}}
 + \frac{l_{1}\ln\left(D_{2}(\nu_{1,0}+l_{1})\right)}{\nu_{1,0}+l_{1}}.
\end{align*}
Next, we only prove the convergence of $\sum_{i=1}^{\infty}\frac{l_{1}\ln\left(D_{2}(\nu_{1,i}+l_{1})\right)}{\nu_{1,i}+l_{1}}$.
By Lemma \ref{lemma6.1} and Lemma \ref{lemma6.2}, we can obtain that
\begin{align*}
   D_{2}(\nu_{1,k}+l_{1})
&=  D_{2}\frac{a^{k+1}-1}{a-1}\nu_{1,0}+D_{2}l_{1}
=  \frac{D_{2}\nu_{1,0}a}{a-1}a^{k}-\frac{D_{2}\nu_{1,0}}{a-1}+D_{2}l_{1}
\leq  \frac{D_{2}\nu_{1,0}a}{a-1}a^{k}+D_{2}l_{1}
   \nonumber\\
&
\leq  \max\left\{\left(\frac{D_{2}\nu_{1,0}a}{a-1}\right)^{k},\left(D_{2}l_{1}\right)^{k}\right\}a^{k}
\end{align*}
Thus, we have
\begin{align}\label{6.1.2}
     \ln(D_{2}(\nu_{1,k}+l_{1}))
\leq k\ln\left( \max\left\{\frac{D_{2}\nu_{1,0}a}{a-1},D_{2}l_{1}\right\}a\right).
\end{align}
Besides, by Lemma \ref{lemma6.1} and Lemma \ref{lemma6.2}, we also have
\begin{align}\label{6.1.3}
   \nu_{1,k}+l_{1}
&=  \frac{a^{k+1}-1}{a-1}\nu_{1,0}+l_{1}
=   (a^{k}+\cdots+a+1)\nu_{1,0}+l_{1}
\geq a^{k}\nu_{1,0}.
\end{align}
So, by \eqref{6.1.2} and \eqref{6.1.3}, we have
\begin{align*}
     \sum_{k=1}^{n}\frac{l_{1}\ln\left(D_{2}(\nu_{1,k}+l_{1})\right)}{\nu_{1,k}+l_{1}}
\leq \sum_{k=1}^{n}\frac{
      kl_{1}\ln\left( \max\left\{\frac{D_{2}\nu_{1,0}a}{a-1},D_{2}l_{1}\right\}a\right)}{a^{k}\nu_{1,0}}
=    \frac{
      l_{1}\ln\left( \max\left\{\frac{D_{2}\nu_{1,0}a}{a-1},D_{2}l_{1}\right\}a\right)}{\nu_{1,0}}
      \sum_{k=1}^{n}\frac{k}{a^{k}}.
\end{align*}
Thus,
$\sum_{i=1}^{\infty}\frac{l_{1}\ln\left(D_{2}(\nu_{1,i}+l_{1})\right)}{\nu_{1,i}+l_{1}}$
is convergent by the convergence of
$\sum_{k=1}^{\infty}\frac{k}{a^{k}}$. \qed

\begin{lemma}\label{lemma6.4}
 Let $\beta_{k}=\Pi_{i=0}^{k}\frac{1+\nu_{1,i}}{\nu_{1,i}+l_{1}}$.
Then  $0<\beta=\Pi_{i=0}^{\infty}\frac{1+\nu_{1,i}}{\nu_{1,i}+l_{1}}<1$, as $k\rightarrow\infty$.
\end{lemma}
\par
\noindent
{\bf Proof.}
First, we can get the following inequality
\begin{align}\label{6.1.4}
\nu_{1,0}+l_{1}-a>0.
\end{align}
In fact, by
$a=\frac{l_{1}^{\ast}}{ql_{1}}$ and $\nu_{1,0}=\frac{l_{1}^{\ast}-q}{q}$, we have
\begin{align*}
\nu_{1,0}+l_{1}-a
=\frac{l_{1}^{\ast}-q}{q}+l_{1}-\frac{l_{1}^{\ast}}{ql_{1}}
=\frac{ql_{1}^{2}+l_{1}^{\ast}l_{1}-ql_{1}-l_{1}^{\ast}}{ql_{1}}
=\frac{ql_{1}(l_{1}-1)+l_{1}^{\ast}(l_{1}-1)}{ql_{1}}
>0
\end{align*}
because of $l_{1}>1$.
Moreover, by  $\nu_{1,k}=\frac{1}{q}\left[\frac{(\nu_{1,k-1}+l_{1})N}{N-l_{1}}-q\right]$, $k=1,2,\cdots$, we can get
$\nu_{1,k}+1=a(\nu_{1,k-1}+l_{1})$, $k=1,2,\cdots$.
Thus, by (\ref{6.1.3}), we have
\begin{align*}
 \lim_{k\rightarrow \infty}\beta_{k}
&=\lim_{k\rightarrow \infty}\Pi_{i=0}^{k}\frac{1+\nu_{1,i}}{\nu_{1,i}+l_{1}}
=\lim_{k\rightarrow \infty}
 \left(\frac{1+\nu_{1,0}}{\nu_{1,0}+l_{1}}\Pi_{i=1}^{k}\frac{1+\nu_{1,i}}{\nu_{1,i}+l_{1}}\right)
=\lim_{k\rightarrow \infty}
 \left(\frac{1+\nu_{1,0}}{\nu_{1,0}+l_{1}}\Pi_{i=1}^{k}\frac{a(l_{1}+\nu_{1,i-1})}{\nu_{1,i}+l_{1}}\right)\nonumber\\
&
=\lim_{k\rightarrow \infty}
 \left(\frac{1+\nu_{1,0}}{\nu_{1,0}+l_{1}}\frac{a(l_{1}+\nu_{1,0})}{\nu_{1,k}+l_{1}}\right)
=\lim_{k\rightarrow \infty}\frac{a(1+\nu_{1,0})}{\nu_{1,k}+l_{1}}
=\lim_{k\rightarrow \infty}\frac{a(1+\nu_{1,0})}{(a^{k}+\cdots+a+1)\nu_{1,0}+l_{1}}.
\end{align*}
It is easy to see that
$0<\frac{a(1+\nu_{1,0})}{(a^{k}+\cdots+a+1)\nu_{1,0}+l_{1}}\leq\frac{a(1+\nu_{1,0})}{(a+1)\nu_{1,0}+l_{1}}<1$ by $a>1$, $\nu_{1,0}>0$, $l_{1}>1$ and \eqref{6.1.4}.
So, we have
$0<\lim_{k\rightarrow \infty}\beta_{k}<1$, i.e.,
$0<\beta=\Pi_{i=0}^{\infty}\frac{1+\nu_{1,i}}{\nu_{1,i}+l_{1}}<1$. \qed
\vskip3mm
 \noindent
\noindent{\bf Funding information}

\noindent
This work is supported by Yunnan Fundamental Research Projects (grant No: 202301AT070465), Xingdian Talent Support Program for Young Talents of Yunnan Province, and Guangdong Basic and Applied Basic Research Foundation (No: 2020A1515110706).

\bibliographystyle{amsplain}

 \end{document}